\documentclass[11pt]{amsart}
\usepackage{amscd,amssymb,amsmath,latexsym,enumerate}
\usepackage{bbm}
\usepackage[mathscr]{euscript}
\usepackage{epsfig}
\usepackage{color}

\newtheorem{theo}{Theorem}
\newtheorem{defini}{Definition}
\newtheorem{proposi}{Proposition}
\newtheorem{lemma}{Lemma}
\newtheorem{coro}{Corollary}
\newtheorem{rem}{Remark}
\newtheorem{exam}{Example}

\newtheorem{prob}{Problem}

\newcommand{\Cc}{{\mathcal C}}

\newcommand{\Ee}{{\mathcal E}}
\newcommand{\Ff}{{\mathcal F}}

\newcommand{\Pp}{{\mathcal P}}

\newcommand{\Uu}{{\mathcal U}}

\newcommand{\CM}{{\mathbb C}}

\newcommand{\NM}{{\mathbb N}}

\newcommand{\RM}{{\mathbb R}}

\newcommand{\ZM}{{\mathbb Z}}

\newcommand{\AG}{{\mathfrak A}}

\newcommand{\as}{{\mathscr A}}
\newcommand{\bs}{{\mathscr B}}
\newcommand{\cs}{{\mathscr C}}

\newcommand{\es}{{\mathscr E}}  

\newcommand{\gs}{{\mathscr G}}

\newcommand{\js}{{\mathscr J}}
\newcommand{\ks}{{\mathscr K}}

\newcommand{\os}{{\mathscr O}}

\newcommand{\us}{{\mathscr U}}
\newcommand{\vs}{{\mathscr V}}
\newcommand{\ws}{{\mathscr W}}

\newcommand{\te}{{\tilde{e}}}

\newcommand{\tu}{{\tilde{u}}}
\newcommand{\tv}{{\tilde{v}}}

\newcommand{\Cs}{$C^{\ast}$-algebra }           %%C*-algebra w. end space
\newcommand{\CsS}{$C^{\ast}$-algebras}          %%C*-algebras
\newcommand{\Orb}{\textit{Orb}}                 %%Orbit
\newcommand{\pea}{\Pp\Ee}			%%Pattern equiv fctions
\newcommand{\tra}{\mbox{\sc t}}                 %%translation-shift
\newcommand{\tri}{\mbox{\sc \tiny t}}           %%shift index
\newcommand{\CB}{\overline{\mathrm B}}          %%Closed Ball
\newcommand{\Inv}{\js}                          %%closed invariant sets
\newcommand{\Dic}{\mathfrak W}               	%%space of dictionaries
\newcommand{\dic}{\ws}               		%%space of dictionaries
\newcommand{\scf}{\wp}	  			%%subword complexity
\newcommand{\csp}{\as^\ZM}	  		%%configuration space
\newcommand{\GAP}{\textbf{G}}	  		%%Gap graphs
\newcommand{\nbr}{N_b}	  			%%# branching vertices
\begin{document}

\title[Spectral continuity for aperiodic quantum systems II]{Spectral continuity for aperiodic quantum systems II.\\
        Periodic approximations in 1D}
\thanks{Work supported in part by NSF Grant No. 0901514 and DMS-1160962 and SFB 878, M\"unster. G.D.'s research is supported by the  FONDECYT grant \emph{Iniciaci\'{o}n en Investigaci\'{o}n 2015} - $\text{N}^{\text{o}}$ 11150143.}

\author{Siegfried Beckus, Jean Bellissard, Giuseppe De Nittis}

\address{Department of Mathematics\\
Technion - Israel Institute of Technology\\
Haifa, Israel}
\email{beckus.siegf@technion.ac.il}

\address{Georgia Institute of Technology\\
School of Mathematics\\
Atlanta GA, USA}
\address{Westf\"alische Wilhelms-Universit\"at\\
Fachbereich 10 Mathematik und Informatik\\
M\"unster, Germany}
\email{jeanbel@math.gatech.edu}

\address{
Facultad de Matem\'aticas \& Instituto de F\'{\i}sica\\
Pontificia Universidad Cat\'olica\\
Santiago de Chile, Chile}
\email{gidenittis@mat.uc.cl}

\begin{abstract}
The existence and construction of periodic approximations with convergent spectra is crucial in solid state physics for the spectral study of corresponding Schr\"odinger operators. In a forthcoming work \cite{BBdN17} this task was boiled down to the existence and construction of periodic approximations of the underlying dynamical systems in the Hausdorff topology. As a result the one-dimensional systems admitting such approximations are completely classified in the present work. In addition explicit constructions are provided for dynamical systems defined by primitive substitutions covering all studied examples such as the Fibonacci sequence or the Golay-Rudin-Shapiro sequence. One main tool is the description of the Hausdorff topology by the local pattern topology on the dictionaries as well as the GAP-graphs describing the local structure. The connection of branching vertices in the GAP-graphs and defects is discussed.
\end{abstract}

%%%%%%%%%%%%%%%%%%%%%%%%%%%%%%%%%%%%%%%%%%%%%%%%%%%%%%%%%%%%%%%%%%%%

\maketitle
\tableofcontents

%%%%%%%%%%%%%%%%%%%%%%%%%%%%%%%%%%%%%%%%%%%%%%%%%%%%%%%%%%%%%%%%%%%%%%%%%

\section{Introduction and main results}
\label{1D.sect-MainResult}

\noindent The aim of the present paper is to characterize when the spectrum of a Hamiltonian in dimension one can be approximated by the spectrum of periodic Hamiltonians in the Hausdorff metric. In this work, only the spectrum as a set will be investigated. In particular, no result on the nature of the spectral measures will be offered. However, the existence and construction of such approximations has been addressed in solid state physics \cite{BS91,Pr13} and a general mathematical method to deal with such problems was missing so far. Even in dimension one, where most efforts from Mathematicians were provided, the present paper is bringing more results on this issue than the ones found in literature. As it turns out the one-dimensional case can be reduced to operators defined by a symbolic dynamical system over a finite alphabet (Proposition~\ref{1D.prop-EncDel}). Dealing with approximations requires to face the problem of creating defects in the spectrum \cite{OK85,BIT91}. As promised in the previous paper \cite[Section~1.5.4]{BBdN17}, an algorithmic method is provided in the present article to construct periodic approximations immune to creating defects (Theorem~\ref{1D.th-ExPerAppr} and Proposition~\ref{1D.rem-PrimSub}). In addition, a complete classification of systems admitting periodic approximations is given (Theorem~\ref{1D.th-ExPerAppr}). As it turns out all relevant examples found in the literature admit periodic approximations (see \cite{BBdN18} for a list of examples). Periodic approximations have spectrum computable through the Floquet-Bloch theory and many software are now available to implement such calculation numerically. For the sake of length, only the Fibonacci sequence and the Golay-Rudin-Shapiro sequence are considered in this paper as a toy models. A book of examples is discussed separately in \cite{BBdN18}.

\vspace{.1cm}

\noindent The spectral convergence is a byproduct of the construction of a continuous field of operators \cite{DD63,Di69}. In fact, having a continuous field of operators is enough to imply the continuity of the spectra as sets and the continuity of the gap edges \cite{BB16}. Due to the fact that the approximating operators are periodic they have an absolutely continuous spectra and one might ask about the spectral type of the limit operator. As a matter of fact, the type of spectral convergence associated with a continuous field of operators is too weak for preserving finer information like the  pure point, singular continuous or absolutely continuous, nature of the spectral measures. Thus, some extra tool is needed to reconstruct the spectral measure. This question is not investigated here.

%%%%%%%%%%%%%%%%%%%%%%%%%%%%%%%%%%%%%%%%%%%%%%%%%%%%%%%%%%%%%%%%%%%%%%%%%

\subsection{Framework and results}
\label{1D.ssect-resu}

Solids are mathematically modeled by so called {\em Delone sets} \cite{Be15}. In addition, it is convenient to assume {\em finite local complexity} meaning that at most finitely many different patches of atomic configuration can appear in a given ball. As shown in Proposition~\ref{1D.prop-EncDel}, this assures that any one-dimensional Delone set of finite local complexity can be encoded by a two-sided infinite word associated with a finite alphabet $\as$. Specifically, it suffices to consider the configuration space $\csp:=\prod_{n\in\ZM}\as=\{\xi:\ZM\to\as\}$. The latter set can be equipped with the product topology and the group $\ZM$ acts continuously on $\csp$ by translations induced by the shift $\tra^m(\xi)(n):=\xi(n-m)$. Thus, $(\csp,\ZM,\tra)$ is a topological dynamical system \cite{GH55}. A subset $\Xi\subseteq\csp$ is called {\em invariant} if $\tra^m(\Xi)\subseteq\Xi$ for all $m\in\ZM$. 
A closed, $\ZM$-invariant subset of $\as^\ZM$ is called a {\em subshift} and the set of all non-empty subshifts is denoted by $\Inv$. The set $\Inv$ equipped with the Hausdorff topology gets the structure of a second countable compact Hausdorff space \cite{Bec16,BBdN17}.

\vspace{.1cm}

\noindent A special class of discrete Schr\"odinger operators, called {\em finite range Hamiltonians}, is built as follows. Let $\ell^2(\ZM)$ be the Hilbert space on which this operator acts. 
Then $\ZM$ is represented by its left regular representation defined by

$$U_m\psi(n)=\psi(n-m)\,,
   \hspace{2cm}
    \psi\in\ell^2(\ZM)\,,\, n,m\in\ZM\,.
$$

\noindent Let $K\subseteq \ZM\setminus\{0\}$ be a finite subset. Let $v:\csp\to\CM$ and $t_k:\csp\to\CM$ for $k\in K$ be continuous functions. For $\xi\in\csp$, the {\em Hamiltonian} $H_\xi:\ell^2(\ZM)\to\ell^2(\ZM)$ is defined by

%%%%%%%%%%%%%%%%
\begin{equation}
\label{1D.eq-schdis}
(H_\xi\psi)(n) \; := \; 
  \sum\limits_{k\in K} 
   t_k(\tra^{-n}(\xi)) \cdot \psi(n-k)
    + v(\tra^{-n}(\xi)) \cdot \psi(n)\,,
\end{equation}
%%%%%%%%%%%%%%%%

\noindent where $\psi\in\ell^2(\ZM)$ and $n\in\ZM$. By introducing the hopping functions $\hat{t}_{\xi,k}(n):=t_k(\tra^{-n}(\xi))$ and the potential $\hat{v}_\xi(n):=v(\tra^{-n}(\xi))$, the Hamiltonian \eqref{1D.eq-schdis} can be shortly expressed by
%%%%%%%%%%%%%%%%
\begin{equation}
\label{1D.eq-schdis-bis}
H_\xi \; = \; 
  \sum\limits_{k\in K} 
   \hat{t}_{\xi,k}\; U_k 
    + \hat{v}_\xi\,.
\end{equation}
%%%%%%%%%%%%%%%%
The operator $H_\xi$ is linear and bounded. In order to be self-adjoint the following are required:

(R1) $v$ is real valued,

(R2) $K$ is invariant by $k\mapsto -k$,

(R3) the functions $t_k$ satisfies $t_{-k}(\xi)= \overline{t_k(\tra^{-k}(\xi))}$.

\noindent The family of operators $H:=(H_\xi)_{\xi\in\csp}$ is obviously strongly continuous with respect to the variable $\xi\in\csp$. In addition, it is $\ZM$-covariant, namely

$$U_m H_\xi U_m^{-1}= H_{\tra^m(\xi)}\,,
   \hspace{2cm}
    \xi\in\csp\,,\; m\in\ZM\,.
$$

\noindent In particular, $\sigma(H_\xi)=\sigma(H_{\tra^m(\xi)})$ for all $m\in\ZM$. From strong continuity, it follows that if $\eta$ belongs to the closure of the $\ZM$-orbit $\Orb(\xi)=\{\tra^n\xi\,;\, n\in\ZM\}$ of $\xi$ (this closure is , therefore, a subshift), then $\sigma(H_\eta)\subseteq \sigma(H_\xi)$. 

\vspace{.1cm}

\noindent If $\Xi\subseteq\csp$ is closed and $\ZM$-invariant, then $(\Xi,\ZM,\tra)$ is also a topological dynamical system. In what follows $H_\Xi$ will denote the family $H_\Xi:=(H_\xi)_{\xi\in\Xi}$ and its spectrum $\sigma(H_\Xi)$ will be defined as the union $\sigma(H_\Xi)=\bigcup_{\xi\in\Xi}\sigma(H_\xi)$. It is worth reminding though, that such definition was justified in previous works by its interpretation in terms of \CsS. Namely let $\AG$ denote the \Cs generated by fields of finite range Hamiltonians $H=(H_\xi)_{\xi\in\as^\ZM}$. If the real line $\RM$ is equipped with its canonical Euclidean metric, let $\ks(\RM)$ denotes the set of compact subsets of $\RM$, equipped with the corresponding Hausdorff metric. The following result is proven in \cite{BBdN17}, Theorem~2.

%%%%%%%%%%%%%%%%%%
\begin{theo}[\cite{BBdN17}]
\label{1D.th-ContSpectrGenSchrDyn}
\noindent Let $H=(H_\xi)_{\xi\in\as^\ZM}$ be a field of self-adjoint operators defining an element of the \Cs $\AG$. For each $\Xi\in \Inv$, let $H_\Xi$ denotes its restriction to $\Xi$. Then the map 

$$\Sigma_H:\Inv\to\ks(\RM)\,,
	\quad 
   \Xi\mapsto \sigma(H_\Xi)\,,
$$

\noindent is continuous in the corresponding Hausdorff topologies.
\end{theo}
%%%%%%%%%%%%%%%%%%

\noindent In particular, the set $K$ defined in Equation~(\ref{1D.eq-schdis}), may be taken infinite provided the hopping functions $\hat{t}_{\xi,k}$ decay sufficiently quickly at infinity in $k$ to define an element of $\AG$.  

\vspace{.1cm}

\noindent The second step consists in answering the following question: When is it possible to approximate the Hamiltonian $H_\Xi$ for some $\Xi\in\Inv$ by periodic Hamiltonians $(H_{\Xi_k})$? By Theorem~\ref{1D.th-ContSpectrGenSchrDyn}, this boils down to the existence of periodic approximations $(\Xi_k)_k$ of the subshift $\Xi$ in the Hausdorff topology on $\Inv$. It is worth reminding that an element $\xi\in\as^\ZM$ is periodic whenever there is a natural integer $n\in\NM$ such that $\tra^n\xi=\xi$. In this case the minimum such integer, $q$ is called the period of $\xi$ and it is sufficient to fix a subword of length $q$ expressed in the alphabet $\as$ to define the sequence $\xi$ entirely. This shows in turn that $\Orb(\xi)=\Xi$ is finite, thus closed, and minimal. It will be called a {\em periodic subshift}. It is easy to check that if, in turn, $\Orb(\xi)=\Xi$ is finite, then it is minimal and periodic (more generally, see \cite{BBdN17} and Definition~11 and Proposition~3 in Section~3.1). In this spirit, a subshift $\Xi\in\Inv$ will be called {\em periodically approximable}, if there is a sequence $\Xi_k\in\Inv$ of periodic subshifts such that $\Xi_k\to\Xi$ in the Hausdorff topology of $\Inv$. It is worth noticing that not any subshift is periodically approximable (see Example~\ref{1D.exam-wall}).

\vspace{.1cm}

\noindent Thanks to the Floquet-Bloch theory, every Hamiltonian $H_\Xi$ associated with a $\Xi$ which is periodically approximable, admits a sequence of algorithmic computable spectra. Within this work, we are going to characterize the property of periodically approximable subshifts as follows: to every subshift $\Xi$, the so called GAP-graphs $\GAP=(\gs_k)_{k\in\NM}$ are defined (See Section~\ref{1D.sect-GAP}) where the vertices are words of length $k$ appearing in $\Xi$. An oriented edge from the vertices $u$ to $v$ is defined if there is a word $w$ of length $k+1$ such that the prefix of length $k$, denoted by $\partial_0 w$, of $w$ is $u$ and the suffix of length $k$, denoted by $\partial_1 w$, of $w$ is $v$. Using the left to right order in words, $\gs_k$, equipped with the boundary maps $\partial_0,\partial_1$, becomes an oriented graph. An oriented graph is called {\em strongly connected} if for any pair {of} vertices $(u,v)$ there are two oriented paths $\gamma_1$ and $\gamma_2$  such that $\gamma_1$ joins the vertices $u$ to $v$ and $\gamma_2$ joins the vertices $v$ to $u$, (see Definition~\ref{1D.def-StrongConnect}). 

%%%%%%%%%%%%%%%%
\begin{theo}
\label{1D.th-ExPerAppr} 
Let $\as$ be a finite alphabet. A subshift $\Xi\in\Inv$ is periodically approximable if and only if 
there is a subsequence of GAP-graphs $(\gs_{k_l})_{l\in\NM}$, with $k_l\to \infty$, that are strongly connected. 
\end{theo}
%%%%%%%%%%%%%%%%

\noindent More precisely, if $\Xi\in\Inv$ is periodically approximable, there exists a sequence of periodic elements $\eta_l\in\csp$, $l\in\NM$, such that the associated sequence of subshifts $(\Xi_l)_l$ given by $\Xi_l:=\Orb(\eta_l)$ converges to $\Xi$ in $\Inv$. Here, every periodic element $\eta_l$ is obtained by choosing a global path in the corresponding GAP-graph $\gs_{k_l}$ following the prescription given in Definition~\ref{1D.def-AssPerWor}. 

\noindent A consequence of Theorem \ref{1D.th-ExPerAppr} is the following

%%%%%%%%%%%%%%%%
\begin{coro}
\label{1D.cor-MinAllPath}
Every minimal $\Xi\in\Inv$ is periodically approximable.
\end{coro}
%%%%%%%%%%%%%%%%

\noindent A very practical tool to describe the Hausdorff topology, essential to prove Theorem~\ref{1D.th-ExPerAppr}, is what will be called here the {\em local pattern topology}. Namely two sequences $\xi,\eta$ in $\as^\ZM$ are Hausdorff close if and only if there is an integer $k$ such that any word of length $k$ in one is a word of the other. This defines a topology on the set of dictionary $\Dic$ (Definition~\ref{1D.def-dict}). Then, \cite{Bec16}

%%%%%%%%%%%%%%%%
\begin{theo}[\cite{Bec16}]
\label{1D.th-DictHomeoSubs}
Let $\as$ be a finite alphabet. The topological spaces $\Inv$ and $\Dic$ are homeomorphic.
\end{theo}
%%%%%%%%%%%%%%%

%%%%%%%%%%%%%%%%
\begin{rem}
\label{1D.rem-Gshift}
{\em The notion of subshift and dictonary (see \cite{Bec16}) extends naturally to $\as^G$ where $G$ is a discrete countable group. The existence of a bijection between $\js$ and $\Dic$ was known before for $G =\ZM^d$ \cite{LM95}. The topological aspect of an homeomorphism was first
proven in \cite{Bec16} for general countable groups G.
}%%
\hfill $\Box$
\end{rem}
%%%%%%%%%%%%%%%%

\noindent Substitutions played a crucial role in the recent past to define various subshifts. Some standard references are \cite{Qu87,Fog02,Qu10}. Given an alphabet $\as$, let $\as^\ast$ be the set of all {\em words} of finite lengths with letters in $\as$. This set has an associative product obtained by concatenation of two words and the empty word \o\, plays the role of a neutral element. A {\em substitution} is a map $S:\as^\ast\to\as^\ast$ such that $S(vw)=S(v)S(w)$ (homomorphism). In particular any substitution is defined once the words $\{S(a)\,:\, a\in\as\}$ have been defined. It is customary to restrict the choice of a substitution to the following class: a substitution is called {\em primitive} if there is $N\in\NM$ such that for any pair $a,b$ of letters, $a$ is a letter in the word $S^N(b)$. In this case a dictionary $\dic(S)$ is defined as the set of all subwords of $S^n(a)$ for some $n\in\NM$ and $a\in\as$ (see Definition~\ref{1D.def-dict}). This dictionary, in turn, defines a unique subshift $\Xi(S)$, thanks to Theorem~\ref{1D.th-DictHomeoSubs}. If $S$ is primitive, then $\Xi(S)$ is minimal \cite{Fog02}. Moreover,

%%%%%%%%%%%%%%%%
\begin{proposi}
\label{1D.prop-PrimSub}
Let $S$ be a primitive substitution and let $\Xi\in\Inv$ be the associated  subshift. Then $\Xi$ is periodically approximable. Furthermore, if $\gamma$ is a closed path in the GAP-graph $\gs_l(\Xi)$ for some $l\in\NM$ with associated periodic word $\eta:=\eta(\gamma)\in\csp$, then $\Xi_k:=\Orb(S^k(\eta))$ defines a sequence of periodic subshifts converging to $\Xi$.
\end{proposi}
%%%%%%%%%%%%%%%%

%%%%%%%%%%%%%%%%%%%%%%%%%%%%%%%%%%%%%%%%%%%%%%%%%%%%%%%%%%%%%%%%%%%%%%%%%

\subsection{A short historical review}
\label{1D.ssect-hist}

\noindent In 1984, Shechtman et.al \cite{SBGC84} discovered an AlMn alloy with sharp diffraction pattern and a ten fold symmetry incompatible with periodic structures in 3D. Such alloys are called quasicrystals nowadays being non-periodic but with long range order. Since then, a large number of quasicrystals alloys were discovered inspiring mathematicians and physicists over the last three decades. In order to study the electronics properties of such materials, a tight-binding representation of the corresponding Schr\"odinger operator is considered, leading to a discrete Schr\"odinger operator of the type proposed in the previous section. From the mathematical point of view, mainly one-dimensional systems were studied so far. The reader is invited to a much more detailed discussion in \cite{Be86,AG95,DLQ14,Be15,DEG15} and references therein.

\vspace{.1cm}

\noindent One of the first example studied in the literature was the {\em Fibonacci Hamiltonian} {which is} a special example of a class of one-dimensional systems over a two-letter alphabet, called {\em Kohmoto model} \cite{KKT83,OPRSS83,KO84,OK85}. In \cite{KO84}, the transfer matrix method and the trace map were used to numerically conclude that the Fibonacci model has a Cantor spectrum. A mathematical breakthrough was obtained by Casdagli \cite{Ca86} and S\"ut\H{o} \cite{Su87,Su89} proving that the spectrum of the Fibonacci Hamiltonian is a Cantor set of Lebesgue measure zero and the spectral measure is purely singular continuous. In addition to the previously mentioned techniques, a Gordon-type argument and explicit periodic approximations were used in an optimal way in \cite{Ca86,Su95}. This methods were tremendously pushed forward over the last decades and extended to a larger class of one-dimensional systems \cite{BIST89,Be90,BBG91,BG93,DL03,DGY16}.

\vspace{.1cm}

\noindent In \cite{OK85}, a numerical study of the Kohmoto model as a function of $\alpha$ (the occurence frequency ratio for one letter) was provided by using rational values of $\alpha$. % (which in turn corresponds to periodic approximation). 
The spectrum of the Kohmoto model was rigorously studied in \cite{BIST89,DL03}. As numerically shown in \cite{OK85}, approximating a rational slope by rational approximation leads to the the creation of a localized defect in the periodic chain. A mathematical proof was provided in \cite{BIT91} showing also that rational approximations of an irrational slope do not create a defect. It is worth mentioning that rational values of the  slope $\alpha$ define periodic Hamiltonians.

\vspace{.1cm}

\noindent The work by Kotani \cite{KO89} provides also a seminal contribution in analyzing the spectrum. More precisely, if the subshift is strictly ergodic and satisfies the Boshernitzan condition then the spectrum equals  the set where the Lyapunov exponent vanishes \cite{Le02,DL06}. Thanks to \cite{KO89}, this leads to a Cantor spectrum of zero Lebesgue measure. At that time, this result covered all known examples including knew models. Later \cite{LQ11} provided subshifts (defined by Toeplitz sequences) that do not satisfy the Boshernitzan condition but admit spectrum of zero Lebesgue measure.

\vspace{.1cm}

\noindent Another big class of one-dimensional quasicrystals is given by primitive substitutions \cite{Fog02}. For such systems purely singular continuous Cantor spectrum of Lebesgue measure zero was shown \cite{Be90,BBG91,BG93,HKS95,Da00b} by exploiting the trace map formalism. This includes examples such as Thue-Morse sequence \cite{Be90,LQY16} and the Period-doubling sequence \cite{Da98,Da00}. In many ways this method corresponds to a periodic approximation, as will be seen in Section~\ref{1D.sect-GAP}. 

\vspace{.1cm}

\noindent Using periodic approximations to compute the spectrum was also promoted in \cite{BS91,Su95} as has been done in \cite{Pr13}. This is because it converges exponentially fast in the period of the approximation as shown in \cite{Pr13}. It is worth noticing that periodic approximations are also important to estimate fractal dimensions of the spectrum, see e.g. \cite{DEGT08,LQ15}.

\vspace{.1cm}

\noindent As discussed before, the methods of transfer matrices and trace maps are powerful tools. On the other hand, these techniques are limited to Hamiltonians with nearest neighbor interaction and to one-dimensional systems or models that can be decomposed in one-dimensional systems \cite{DGS15}. But higher dimensional systems like the Penrose tiling or the Octagonal lattice cannot be treated with this methods \cite{BS91}. As a common methodology of the works in the past decades suitable periodic approximations were used to extract spectral informations of Schr\"odinger operators associated with quasicrystals. The main goal of \cite{Bec16,BBdN17} was the construction of a general theory which allows to characterize the convergence of the spectra for much larger class of Hamiltonians by looking at the convergence of the underlying structures, see also Theorem~\ref{1D.th-ContSpectrGenSchrDyn}. The formalism used in \cite{Bec16,BBdN17} makes extensive use of the theory of groupoids and \CsS. The whole construction is based on the use of the tautological groupoid. The aim is to investigate this approach further. In view of this, the present paper provides a detailed description of this theory for the one-dimensional case. In addition this convergence of the underlying structures implies the weak-$\ast$ convergence of measured quantities such as the density of state measure under suitable assumptions \cite{BP17}.

\vspace{.1cm}

\noindent Finally it is worth mentioning the recent paper \cite{KP17} which focuses on the study of the stability of edge states produced by cutting half of the system in Sturmian models. The analysis of these states is performed by adapting the approximation scheme developed in \cite{Bec16,BBdN17} to half-space systems. This is, in particular, an evidence that the method described here is suitable for generalizations and applications to various problems of interest in condensed matter. 

%%%%%%%%%%%%%%%%%%%%%%%%%%%%%%%%%%%%%%%%%%%%%%%%%%%%%%%%%%%%%%%%%%%%%%%%%

\subsection{Organization of the paper}
\label{1D.ssect-Orga}

\noindent The encoding of Delone sets of finite local complexity as two-sided infinite sequences over an finite alphabet is discussed in Subsection~\ref{1D.ssect-DelFLC}. Afterwards, the notions of dictionaries (Subsection~\ref{1D.ssect-Dic}) and subshifts (Subsection~\ref{1D.ssect-Subsh}) are introduced as well as the corresponding topologies. Based on this, the proof of Theorem~\ref{1D.th-DictHomeoSubs} and its consequences are provided in Subsection~\ref{1D.ssect-th3}. The fundamental notion of GAP-graphs is introduced in Section~\ref{1D.sect-GAP} after recalling basic facts and notions of graphs. In Subsection~\ref{1D.ssect-ConstrGAP}, an algorithmic recipe is provided to construct the sequence of GAP-graphs. Fundamental properties of these graphs as well as the proof of Theorem~\ref{1D.th-ExPerAppr} and its consequences are presented in Subsection~\ref{1D.ssect-PropGAP}. Subsection~\ref{1D.ssect-defect} is devoted to a discussion of the connection between the subword complexity and the branching vertices of the GAP-graphs. In light of this, the relation with defects is discussed there. The toy models are presented in Subsection~\ref{1D.ssect-Fib} and Subsection~\ref{1D.ssect-GRS}. Appendix \ref{1D.sect-PattEqOp} is devoted to the presentation of the  pattern equivariant algebra.

\vspace{.5cm}

\noindent {\bf Acknowledgments: } S.B. wants to express his deep gratitude to Daniel Lenz for his constant and fruitful discussions and support over the last years. In addition S.B. would like to thank Tobias Hartnick for pointing out the work by de~Bruijn \cite{Bru46}. G.D. wants   to thank Frederic Klopp and Fran\c{c}ois Germinet   for their support during the  initial years of the investigation which finally leads to this work.
This work has greatly benefited from the support of various institutions, the Mathematics Department at Technion, Israel, the Facultad de Matem\'aticas at the Pontificia Universidad Cat\'olica, Chile, the Department of Mathematics Westf\"alische Wilhelms-Universit\"at, M\"unster, Germany, Georgia Institute of Technology, USA, the Erwin Schr\"odinger Institute, Vienna and the Research Training Group (1523/2) at the Friedrich-Schiller University in Jena, Germany. This research was supported through the program ``Research in Pairs'' by the Mathematisches Forschungsinstitut Oberwolfach in 2018. The work is part of the NSF grant ``Spectral Properties of Aperiodic Solids'', Grant No. DMS1160962. G.D.'s research is supported by the  FONDECYT grant \emph{Iniciaci\'{o}n en Investigaci\'{o}n 2015} - No. 11150143.

\vspace{.3cm}

%%%%%%%%%%%%%%%%%%%%%%%%%%%%%%%%%%%%%%%%%%%%%%%%%%%%%%%%%%%%%%%%%%%%%%%%%

\section{Encoding theory}
\label{1D.sect-enc}

\noindent In this section Delone sets of finite local complexity are defined. This sets provide the standard model for the mathematical description of quasicrystals and aperiodic solids \cite{Be15}. As it turns out in the one-dimensional case a Delone set of finite local complexity can be encoded combinatorially as a two-sided infinite word.

%%%%%%%%%%%%%%%%%%%%%%%%%%%%%%%%%%%%%%%%%%%%%%%%%%%%%%%%%%%%%%%%%%%%%%%%%

\subsection{Encoding of colored Delone sets with finite local complexity}
\label{1D.ssect-DelFLC}

%%%%%%%%%%%%%%%%%%%%%%%%%%%%%%%%%%%%%%%%%%%%%%%%%%%%%%%%%%%%%%%%%%%%

\noindent Let $D\subset \RM$ be a discrete set. It is called {\em uniformly discrete} whenever there is an $r>0$ such that in any open Euclidean ball of radius $r$ there is at most one point of $D$. It is called {\em relatively dense} if there is $R>0$ such that in any closed Euclidean ball of radius $R$ there is at least one point of $D$. It is called a {\em Delone set} if it is both uniformly discrete and relatively dense. A {\em patch of radius $s>0$} is a finite subset of $\RM$ of the form $p=\CB(0;s)\cap (D-x)$ for some $x\in D$ where $\CB(0;s)$ denotes a closed Euclidean ball of radius $s$ around the origin $0$. A Delone set is called {\em repetitive} whenever, given any patch $p$ of radius $s$ of $D$, and for any $\epsilon >0$, there is $r_p>0$ such that in any closed Euclidean ball of radius $r_p$, there is a point $x\in D$ such that the Hausdorff distance between $p$ and $\CB(0;s)\cap(D-x)$ is less than $\epsilon$. A Delone set $D$ has {\em finite local complexity}, whenever, for any $s>0$, the number of patches of $D$ of radius $s$ is finite. A Delone set of finite local complexity is called {non-periodic} whenever $D+x=D$ for $x\in\RM$ implies $x=0$.

\vspace{.3cm}

\noindent In view of $D$ representing the positions of the atomic nuclei in a solid with several atomic species, it will be assumed that there is a finite set $\cs$ called the {\em color set}. A Delone set $D\subset \RM$ with a map $x\in D\mapsto c(x)\in\cs$ is called a {\em colored Delone set}.

\vspace{.3cm}

\noindent The following assertion provides a description of $D$ as a two-sided infinite word with letters in a finite alphabet.

%%%%%%%%%%%%%%%%
\begin{proposi}
\label{1D.prop-EncDel}
A Delone set $D$ of finite local complexity can be encoded as a two-sided infinite word over a finite alphabet. Moreover, the word length metric is equivalent to the Euclidean metric on a Delone set. 
\end{proposi}
%%%%%%%%%%%%%%%%

\noindent  {\bf Proof: } (i) Let $x_0$ be the minimum $\min\{D\cap[0,\infty)\}$. Define inductively $x_{\pm n}$ as $x_{n+1}=\min\{D\cap(x_n,+\infty)\}$ and $x_{-n}=\max\{D\cap (-\infty,x_{-n+1})\}$ for $n\in\NM$. Then, $D=(x_n)_{n\in\ZM}$ with $x_n<x_{n+1},\; n\in\ZM$. The alphabet associated with the colored Delone set $D$ of  finite local complexity is defined by

$$\as_D
	:=\{ (x_{n+1}-x_n,c(x_n))\;:\; n\in\ZM\}.
$$ 

\noindent Since $D$ is Delone it follows that $2r\leq x_{n+1}-x_n\leq 2R$. Using the finite local complexity and the finiteness of $\cs$ the set $\as_D$ is finite. Each letter $a=(\ell_a,c(a))\in\as_D$ can be seen as an interval of length $[0,\ell_a]$ punctured at $0$ having a color $c(a)$. For $D=(x_n)_{n\in\ZM}$ a two-sided infinite word $\xi_D:\ZM\to\as_D$ is associated where $\xi_D(n):=(x_{n+1}-x_n,c(x_n))$ for $n\in\ZM$.

\vspace{.1cm}

\noindent (ii) Conversely, given a two-sided infinite word $\xi$ over the alphabet $\as$. Then a colored Delone set of finite local complexity can be reconstructed by identifying a letter $a\in\as$ with a length $\ell_a$ and a color $c_a$. By fixing $x_0:=0,\; c(x_0):=c_{\xi(0)}$ and $x_1:=\ell_{\xi(0)}$ such a Delone set could be reconstruct  inductively. 

\vspace{.1cm}

\noindent The previous encoding shows that instead of using the Euclidean metric to measure distances on $D$ it is possible to use the word length $d(x_n,x_m)=|n-m|$. The Delone condition implies

$$2r|n-m|\leq |x_n-x_m|\leq 2R|n-m|\,.
$$

\noindent Hence, the two metrics are equivalent.
\hfill$\Box$

%%%%%%%%%%%%%%%%%%%%%%%%%%%%%%%%%%%%%%%%%%%%%%%%%%%%%%%%%%%%%%%%%%%%%%%%%
\subsection{Dictionaries}
\label{1D.ssect-Dic}

\noindent An \emph{alphabet} is a \emph{finite} set $\as$. A \emph{word} with letters in $\as$ is an element of the Cartesian product $\as^n$. For the sake of convenience a word $u=(a_1,a_2,\cdots,a_n)\in\as^n$ is represented as $u=a_1 a_2\ldots a_n$ with $a_k\in\as$. The number of letters $|u|:=n$ is called its {\em word length}. The \emph{empty word} \o\, has by definition zero length. Given two words $u=a_1\ldots a_n$ and $v=b_1\ldots b_m$ the {\em concatenation} $uv=a_1\ldots a_n b_1\ldots b_m$ is an element of $\as^{n+m}$. Then the concatenation map $(u,v)\mapsto uv$ is associative, but not commutative. For a word $u$ and $j\in\NM$ the $j$-time concatenation of the word $u$ is denoted by $u^j:=uu\ldots u\in\as^{j|u|}$. A word $u$ is called a {\em subword} of a word $v$ if there exist two other words $w_1,w_2$ (possibly the empty words) such that $v=w_1uw_2$. In particular, this implies that $|u|\leq |v|$.

%%%%%%%%%%%%%%%%
\begin{defini}[The Space of Dictionaries]
\label{1D.def-dict}
A non-empty family $\dic$ of finite words made by letters of the alphabet $\as$ is called a {\em dictionary} if:
\begin{itemize}
	\item[(D.1)]  \o\,$\in\dic$; \hfill{\bf \upshape (null-element)}
	
	\item[(D.2)] for all $v\in\dic$ each subword $u$ of $v$ belongs to $\dic$; \hfill{\bf \upshape (heredity)}

	\item[(D.3)] for all $u\in\dic$ there exists $a,b\in\as$  such that $aub\in\dic$. \hfill{\bf \upshape (extensibility)}
\end{itemize}

\noindent The set of all dictionaries $\Dic(\as)$ over the alphabet $\as$ is called the space of dictionaries. 
\end{defini}
%%%%%%%%%%%%%%%%

\noindent Whenever there is no ambiguity, the notation $\Dic$ for the space of dictionaries $\Dic(\as)$ over the alphabet $\as$ will be preferred. The set $\Dic$ becomes a topological space if endowed with the {\em local pattern topology} defined by the following basis of sets

%%%%%%%%%%%%%%%%
\begin{equation}
\label{1D.eq-DicTop}
\vs(n,U)
	:=\left\{ \dic\in\Dic
		\;:\; \dic\cap\as^n=U \right\}
			\,,\qquad n\in\NM\,,\; U\subseteq \as^n\,.
\end{equation}
%%%%%%%%%%%%%%%%

\noindent Note that $\dic\cap\as^n$ is the set of words in $\dic$ with length $n$. By the description of a Delone set as a two-sided infinite word (Subsection~\ref{1D.ssect-DelFLC}) a patch can be expressed as a finite word. This satisfies the notion of local pattern topology as the set $\dic\cap\as^n$ is the set of all local patterns up to a certain size $n$.

%%%%%%%%%%%%%%%%
\begin{proposi}
\label{1D.prop-DAprop}
The family $\bs:=\{\vs(n,U)\;:\; n\in\NM,\; U\subseteq\as^n\}$ defines a basis for a topology. Furthermore, $\Dic(\as)$ is second countable, compact, Hausdorff totally disconnected and metrizable in the local pattern topology.
\end{proposi}
%%%%%%%%%%%%%%%%

\noindent  {\bf Proof: } (i) By inspection the sets in Equation~\eqref{1D.eq-DicTop} define a basis. Due to the finiteness of the alphabet $\as$ the set $\as^n$ is finite. Consequently, the set of subsets of $\as^n$ is finite for any fixed $n\in\NM$. Since the countable union of finite sets is countable the basis defined in  Equation~\eqref{1D.eq-DicTop} is countable. In particular, $\Dic$ is second countable. 

\vspace{.1cm}

\noindent (ii) The space $\Dic$ is Hausdorff. Let  $\dic_1$ and $\dic_2$ be two dictionaries. The relation $\vs(n,\dic_1\cap\as^n)\cap\vs(n,\dic_2\cap\as^n)\neq\emptyset$ holds if and only if ${\dic}_1\cap\as^n=\dic_2\cap\as^n$. This equality holds for all $n\in\NM$ if and only $\dic_1=\dic_2$ implying the Hausdorff property. 

\vspace{.1cm}

\noindent (iii) The topological space $\Dic$ is compact. For indeed, let $(\dic_k)_{k\in\NM}$ be a sequence of dictionaries. Then, a convergent subsequence can be extract as follows: since $\as$ is finite there is an $a\in\as$ such that the set $\{ k\in\NM \;:\; a\in\dic_k\cap\as \}$ is infinite. Thus, there is a subsequence of $(\dic_{k})_{k\in\NM}$ such that all elements of the subsequence contain the letter $a\in\as$. By repeating, if necessarily, with the set $\as\setminus\{a\}$ a subsequence $(\dic_{k_l})_{l\in\NM}$ is constructed such that $\dic_{k_l}\cap\as=\dic_{k_m}\cap\as$ for all $l,m\in\NM$. Define $\sigma:\NM\times\NM\to\NM$ inductively. Set $\sigma(1,l):=k_l$ for $l\in\NM$. A similar argument guarantees the existence of a subsequence $(\dic_{\sigma(n-1,l)_i})_{i\in\NM}$ with $\dic_{\sigma(n-1,l)_i}\cap\as^n=\dic_{\sigma(n-1,l)_j}\cap\as^n,\; i,j\in\NM$ for each fixed $n\in\NM$. Set $\sigma(n,j):= \sigma(n-1,l)_j$ for $j\in\NM$. Using a Cantor type argument one extracts a convergent subsequence of the sequence $(\dic_k)_{k\in\NM}$ in the space $\Dic$. Precisely, the sequence $(\dic_{\sigma(n,1)})_{n\in\NM}$ defines a convergent sequence of dictionaries. The limit dictionary $\dic$ is defined by $\dic\cap\as^n=\dic_{\sigma(n,1)}\cap\as^n$ for $n\in\NM$.

\vspace{.1cm}

\noindent (iv) By (ii) the space $\Dic$ is Hausdorff. Thus, if $\vs(n,U)$ is clopen for all $n\in\NM$ and $U\subseteq\as^n$, the space $\Dic$ is totally disconnected. Consider some $n\in\NM$ and $U\subseteq\as^n$. Let $\dic$ be in the closure $\overline{\vs(n,U)}$, i.e. all neighborhoods of $\dic$ intersect $\vs(n,U)$. The set $\vs(n,\dic\cap\as^n)$ is an open neighborhood of $\dic$. Thus, the intersection $\vs(n,\dic\cap\as^n)\cap\vs(n,U)$ is non-empty implying $U=\dic\cap\as^n$. Hence, $\dic\in\vs(n,U)$ follows and so $\vs(n,U)$ is clopen.

\vspace{.1cm}

\noindent Since $\Dic$ is Hausdorff, compact and second countable it is metrizable \cite{Ur24}. 
\hfill$\Box$

%%%%%%%%%%%%%%%%%%%%%%%%%%%%%%%%%%%%%%%%%%%%%%%%%%%%%%%%%%%%%%%%%%%%%%%%%
\subsection{Subshifts}
\label{1D.ssect-Subsh}

\noindent Let $\as$ be a finite alphabet endowed with the discrete topology. The space 

$$\csp:=
	\prod_{n\in\ZM}\, \as
	= \big\{\xi=\big(\xi(n)\big)_{n\in\ZM} \,:\, \xi(n)\in\as \text{ for all } n\in\ZM \big\}
$$ 

\noindent endowed with the product topology is called the {\em configuration space} over the alphabet $\as$. Due to the finiteness of $\as$ the space $\csp$ defines a compact, second countable, Hausdorff, totally disconnected space \cite{Ty30,Ce37}.

\vspace{.1cm}

\noindent An element in $\csp$, i.e. a {\em two-sided infinite word} $\xi$, is identified with a map $\xi:\ZM\to\as$. For each $\xi\in\csp$ and $l,r\in\ZM$ with $l\leq r$ the restriction $\xi|_{[l,r]}$ denotes the {\em subword} $\xi(l)\xi(l+1)\ldots \xi(r-1)\xi(r)$ of length $r-l+1$. The {\em dictionary $\dic(\xi)$ associated with $\xi\in\csp$} is defined by the set of all finite subwords of $\xi$. Then $\dic(\xi)\in \Dic$ follows immediately according to Definition~\ref{1D.def-dict}.

\vspace{.1cm}

\noindent The homeomorphism $\tra:\csp\to\csp$ defined by

$$(\tra\xi)(j):=\xi(j-1),\qquad j\in\ZM,
$$

\noindent is called the {\em shift}. Then $(\csp,\ZM,\tra)$ defines a topological dynamical system. A subset $\Xi\subseteq\csp$ is called {\em invariant} if $\tra(\Xi)\subseteq\Xi$. 

\vspace{.1cm}

\noindent A closed subset $\Xi$ inherits the topology induced by $\csp$. Thus, $\Xi$ is compact, second countable, Hausdorff, totally disconnected  and metrizable. It follows immediately that the cylinder sets defined by

%%%%%%%%%%%%%%%%
\begin{equation}
\label{1D.eq-clopen}
\os(u,v)
 :=\big\{\xi\in\Xi
  \;:\; \xi|_{[-|u|,|v|-1]}=uv\big\}
   \;,\qquad\quad u,v,uv\in\dic(\Xi)\,,
\end{equation}
%%%%%%%%%%%%%%%%

\noindent define a base of open sets for the product topology on the closed subset $\Xi\subseteq\csp$. Note that each $\os(u,v)$ is a \emph{clopen} set (i.e. open and closed). Indeed let $\pi_{l,r}$ denote the (finite) family of open sets $\os(u',v')$ such that $u',v'$ and $u'v'$ are words in $\dic(\Xi)$ with $|u'|=l\,,\,|v'|=r$. Then  (a) their union over $\pi_{l,r}$ is $\Xi$ itself, (b) the intersection of two such distinct open sets is empty. Hence $\pi_{l,r}$ is a finite partition by open sets, so that the complement of each such set, being a union of open sets, is also open. Hence each of these sets is closed, thus clopen. Consequently $\Xi$ is totally disconnected \footnote{The notation $\Xi$ is chosen since it is the most disconnected letter in the greek alphabet.} (see e.g. \cite{Wa82}).

\vspace{.1cm}

\noindent A non-empty, closed and invariant subset $\Xi\subseteq\csp$ is called a {\em subshift}. Then the pair $(\Xi,\ZM,\tra)$ defines a dynamical system. For $\xi\in\csp$ the associated subshift is defined by $\overline{\Orb(\xi)}$ where the set  $\Orb(\xi):=\{\tra^n\xi\;:\; n\in\ZM \}$ is called the {\em orbit of $\xi$}. A subshift $\Xi$ is called {\em topological transitive} if there is a $\xi\in\Xi$ such that $\Xi=\overline{\Orb(\xi)}$. Furthermore, a subshift $\Xi$ is said to be {\em minimal} if for each $\xi\in\Xi$, the orbit $\Orb(\xi)\subseteq\Xi$ is dense. Clearly, every minimal subshift is topological transitive but the converse is false. If $\ks(\csp)$ denotes the set of closed subsets of $\csp$, the set

$$\Inv 
   :=\big\{ \Xi\subseteq\csp\text{ subshift}\big\} \subseteq \ks(\csp)\,,
$$

\noindent is naturally equipped with the induced Hausdorff topology. For a detailed discussion on the topology, the reader is referred to \cite{Bec16,BBdN17}. Throughout this work, the following description of the Hausdorff topology is used: a base for the Hausdorff topology on the compact subsets $\ks(\csp)$ is defined by the family

$$\Uu(F,\Ff):=
   \big\{\Xi\in \ks(\csp)\,:\, \Xi\cap F=\emptyset\;\,\&\;\, 
    \forall O\in \Ff\,,\, \Xi\cap O\neq \emptyset 
    \big\}\,,
$$

\noindent where $F\subseteq\csp$ is closed and $\Ff$ is a finite family of open subsets. This description goes back to \cite{Vi22,Ch50,Fe62}. Then $\Inv\subseteq\ks(\csp)$ is naturally equipped with the induced Hausdorff topology. Thanks to \cite{BBdN17}, Proposition~1, $\Inv$ is a compact, second countable Hausdorff space if equipped with the Hausdorff topology, (see also \cite{Bec16}, Proposition~3.2.5).

\vspace{.1cm}

\noindent Given a subshift $\Xi\in\Inv$, the {\em associated dictionary} is defined by

$$\dic(\Xi)
   :=\bigcup_{\xi\in\Xi} \dic(\xi)\,.
$$

\noindent Clearly, $\dic(\Xi)$ is an element of $\Dic$ in the sense of Definition~\ref{1D.def-dict}. 

%%%%%%%%%%%%%%%%
\begin{proposi}
\label{1D.prop-DicShi}
Let $\xi\in\csp$. Then $\dic(\xi)=\dic(\tra\xi)$ holds and for every $\eta\in\overline{\Orb(\xi)}$, the inclusion $\dic(\eta)\subseteq\dic(\xi)$ holds. In particular, $\dic\big(\overline{\Orb(\xi)}\big)=\dic(\xi)$ holds for all $\xi\in\csp$.
\end{proposi}
%%%%%%%%%%%%%%%%

\noindent {\bf Proof: } The identity $\dic(\xi)=\dic(\tra\xi)$ is obviously satisfied by definition of a dictionary. Let $(n_k)_k\subseteq\ZM$ be chosen such that $\eta=\lim_{k\to\infty} \tra^{n_k}\xi$. Consider some $v\in\dic(\eta)$. Then there is an $i\in\ZM$ such that $\eta|_{[i,i+|v|-1]}=v$. Furthermore, there is an $l\in\NM$ such that $\tra^{n_l}\xi|_{[i,i+|v|-1]}=\eta|_{[i,i+|v|-1]}$ by the definition of the topology on $\csp$. Thus, $v\in\dic(\xi)$ follows.
\hfill $\Box$

\vspace{.2cm}

\noindent Let $u,v,w$ be finite words over the alphabet $\as$. Writing $u=u_0u_1\cdots u_{q-1}$, throughout this work the notation $\xi=u^\infty\in\csp$ denotes the periodic concatenation of the finite word $u$, i.e., $\xi_{nq+l}=u_l$ for $0\leq l<q-1$. This can be written as $ \ldots u u \cdot u u \ldots$ continued to the left and right with $u$'s. Here the dot $\cdot$ fixes the origin, namely the first letter of the word on the right hand side of the dot is at $0$ by convention. Similarly, $u^\infty w\cdot v^\infty$ denotes the two-sided infinite word $\xi\in\csp$ defined by $\xi=\ldots u u u w\cdot v v v\ldots$ continued to the left with $u$'s and to the right with $v$'s.

\vspace{.1cm}

\noindent A two-sided infinite word $\xi\in\csp$ is called {\em periodic} if there exists an $n\in\NM$ such that $\tra^n\xi=\xi$. The smallest such integer $q\in\NM$ is called the {\em period of $\xi$}. Then $\Orb(\xi)$ is closed and contains exactly $q$ elements. In addition $\dic(\xi)\cap \as^n$ has exactly $q$ words \cite{Fog02}. If $\xi\in\csp$ is not periodic it is called {\em non-periodic}. A subshift $\Xi\in\Inv$ is called {\em periodic} if there is a periodic $\xi\in\Xi$ such that $\Orb(\xi)=\Xi$. Then, $\Xi\in\Inv$ is periodic if and only if $\Xi$ is minimal and finite. Furthermore, $\Xi\in\Inv$ is said to be {\em aperiodic} if there is a non-periodic element $\xi\in\Xi$. A subshift $\Xi\in\Inv$ is called {\em completely aperiodic} if $\Xi$ does not contain any periodic element. The central notion of this work is the following

%%%%%%%%%%%%%%%%
\begin{defini}[Periodically approximable]
\label{1D.def-Perapprox}
A subshift $\Xi\in\Inv$ will be called {\em periodically approximable} if there is a sequence of periodic $(\Xi_k)_k\subseteq\Inv$ converging to $\Xi$ in the Hausdorff topology.
\end{defini}
%%%%%%%%%%%%%%%%

%%%%%%%%%%%%%%%%
\begin{rem}
\label{1D.rem-perapp}
{\em Not all subshifts are periodically approximable as shown by Example~\ref{1D.exam-wall} below.
}%%
\hfill $\Box$
\end{rem}
%%%%%%%%%%%%%%%%

\vspace{.1cm}

%%%%%%%%%%%%%%%%%%%%%%%%%%%%%%%%%%%%%%%%%%%%%%%%%%%%%%%%%%%%%%%%%%%%%%%%%%
 \subsection{Dictionaries and the Hausdorff topology}
 \label{1D.ssect-th3}

  %%%%%%%%%%%%%%%%%%%%%%%%%%%%%%%%%%%%%%%%%%%%%%%%%%%%%%%%%%%%%%%%%%
  \subsubsection{Proof of Theorem~\ref{1D.th-DictHomeoSubs}}
  \label{1D.sssect-TH3}

\noindent Let $\phi:\Inv\to\Dic$ be the map defined by

%%%%%%%%%%%%%%%%
$$
\phi(\Xi)\;:=\;\dic(\Xi)\;:=\;\bigcup_{\xi\in\Xi}\dic(\xi).
$$
%%%%%%%%%%%%%%%%

\noindent By definition, the empty set is not an element of $\Inv$. This map is well-defined, i.e. $\dic(\Xi)$  is a dictionary, according to the Definition~\ref{1D.def-dict}. By \cite{Bec16,BBdN17} and Proposition~\ref{1D.prop-DAprop} both spaces $\Inv$ and $\Dic$ are compact and Hausdorff. Thus, according to \cite[Proposition 8.11]{Que01} it suffices to prove the continuity and the bijectivity of $\phi$ so that $\phi$ defines a homeomorphism. Thus it is sufficient to show that $\phi$ is (i) injective, (ii) surjective and (iii) continuous:

\vspace{.1cm}

\noindent (i) {\bf $\phi$ is injective: }Let $\Xi_1,\Xi_2\in\Inv$ be distinct. Without loss of generality it can be assumed that $\Xi_1\setminus\Xi_2\neq\emptyset$. If $\xi\in\Xi_1\setminus\Xi_2$ there exists an $n_0\in\NM$ such that $\xi|_{[-n_0,n_0]}\not\in \phi(\Xi_2)$ and $\xi|_{[-n_0,n_0]}\in \phi(\Xi_1)$ so that $\phi$ is injective. For otherwise, there would be a sequence $\eta_n\in\Xi_2,\; n\in\NM$ such that $\eta_n|_{[-n,n]}=\xi|_{[-n,n]}$. Hence, by definition of the product topology, $\lim_{n\to\infty}\eta_n=\xi\in\overline{\Xi_2}=\Xi_2$, because $\Xi_2$ is closed, a contradiction.

\vspace{.1cm}

\noindent (ii) {\bf $\phi$ is surjective: }Let $\dic$ be a dictionary and let

$$\Xi(\dic)
	:=\{\xi\in\csp\;:\; \dic(\xi)\subseteq\dic\}
$$

\noindent By construction, $\phi(\Xi(\dic))\subseteq\dic$. Moreover $\Xi(\dic)\neq\emptyset$. For if $w\in\dic$, the axiom (D.3) in Definition~\ref{1D.def-dict} implies the existence of a sequence $(w_n)_{n\geq 0}$ of words in $\dic$ such that $w_{n+1}=a_{n+1}w_nb_{n+1}$ for some letters $a_{n+1}, b_{n+1}\in\as$ and $w_0=w$. In particular the length of $w_n$ is  $|w|+2n$ and $w_n$ can be written as $u_nwv_n$ where $u_n,v_n$ are words in $\dic$ of length $n$. Therefore, for each $n\in\NM$ there is element $\eta_n\in\csp$ such that $\eta_n|_{[-n,|w|+n-1]}=w_n$. By definition of the product topology, $\eta=\lim_{n\to \infty}\eta_n$ exists in $\csp$ and $\dic(\eta)\subseteq \dic$. Hence $\eta\in\Xi(\dic)$. Since $w$ has been chosen arbitrarily in $\dic$, it follows that $\dic\subseteq \phi(\Xi(\dic))$. Hence $\dic=\phi(\Xi(\dic))$.

\vspace{.1cm}

\noindent (iii) {\bf $\phi$ is continuous: }Define the clopen set $\os(u):=\{\xi\in\csp\;:\; \xi_{[1,|u|]}=u\}$ for a finite word $u$ with letters in $\as$. Consider a non-empty open set $\vs(n,U)$ in the local pattern topology where $n\in\NM$ and $U:=\{u_1,\ldots, u_l\}\subseteq\as^n$. To prove the continuity of $\phi$ the preimage $\phi^{-1}(\vs(n,U))$ has to contain a non-empty open set. Let then 

$$F
	:= \underset{j=1}{\overset{l}{\bigcap}} \left(\csp\setminus\os(u_j)\right)
	= \big\{\xi\in\csp\;:\; \forall\; 1\leq j\leq l \text{ such that } \xi_{[1,n]}\neq u_j\big\}
$$ 

\noindent a closed subset of $\csp$. Consider the open subsets $\os_j:=\os(u_j),\; 1\leq j\leq l,$ of $\csp$. It suffices to show the equality

$$\us(F,\Ff)
	:= \big\{ \Xi\in\Inv\;:\; F\cap\Xi=\emptyset  \text{ and }  \Xi\cap\os(u_j)\neq\emptyset \text{ for all } 1\leq j\leq l \big\}
		= \phi^{-1}(\vs(n,U))
$$

\noindent where $\Ff:=\{\os_1,\ldots,\os_l\}$. This is checked as follows: Let $\Xi\in\us(F,\Ff)$. Since $\Xi\cap\os(u_j)\neq\emptyset$ for all $1\leq j\leq l$ it follows that $\{u_1,\ldots, u_l\}\subseteq\dic(\Xi)\cap\as^n$. The converse inclusion $\{u_1,\ldots, u_l\}\supseteq\dic(\Xi)\cap\as^n$ follows by $F\cap\Xi=\emptyset$ and invoking the invariance of $\Xi$. Hence, $\Xi\in \phi^{-1}(\vs(n,U))$ follows implying $\us(F,\Ff)\subseteq\phi^{-1}(\vs(n,U))$. 

\vspace{.1cm}

\noindent The opposite inclusion $\phi^{-1}(\vs(n,U))\subseteq\us(F,\Ff)$ holds. For let $\Xi\in \phi^{-1}(\vs(n,U))$ meaning $\dic(\Xi)\cap\as^n=\{u_1,\ldots, u_l\}$. As $u_j,\; 1\leq j\leq l,$ are elements of $\dic(\Xi)$ it follows that $\Xi\cap\os_j\neq \emptyset$ by the invariance of $\Xi$. Furthermore, the intersection $F\cap\Xi$ is empty showing the desired inclusion: For indeed, otherwise there is a $\xi\in\Xi$ such that $v:=\xi|_{[1,n]}\neq u_j$ for all $1\leq j\leq l$. This leads to $v\in\dic(\Xi)\cap\as^n$ and $v\not\in\{u_1,\ldots, u_l\}$, a contradiction with $\Xi\in \phi^{-1}(\vs(n,U))$.
\hfill$\Box$

  %%%%%%%%%%%%%%%%%%%%%%%%%%%%%%%%%%%%%%%%%%%%%%%%%%%%%%%%%%%%%%%%%%
  \subsubsection{Some Consequences of Theorem~\ref{1D.th-DictHomeoSubs}}
  \label{1D.sssect-consTH3}

%%%%%%%%%%%%%%%%
\begin{coro}
\label{1D.cor-SubSpaTotDisc}
The topological space $\Inv$ of subshifts is totally disconnected.
\end{coro}
%%%%%%%%%%%%%%%%

\noindent {\bf Proof: } 
This follows immediately by Proposition~\ref{1D.prop-DAprop} and Theorem~\ref{1D.th-DictHomeoSubs}. 
\hfill$\Box$

%%%%%%%%%%%%%%%%
\begin{coro}
\label{1D.cor-ConSubDic}
A sequence of subshifts $(\Xi_k)_k\subseteq\Inv$ converges to $\Xi\in\Inv$ if and only if, for every $m\in\NM$, there is a $k_m\in\NM$ such that
$$\dic(\Xi)\cap\as^m
	=\dic(\Xi_k)\cap\as^m
	\,,\qquad k\geq k_m\,.
$$
\end{coro}
%%%%%%%%%%%%%%%%

\noindent {\bf Proof: } 
This is a direct consequence of Theorem~\ref{1D.th-DictHomeoSubs}. 
\hfill$\Box$

%%%%%%%%%%%%%%%%
\begin{exam}
\label{1D.exam-wall}
{\em Not all subshifts are periodically approximable. A good counter example is the following. Let $\as=\{a,b\}$ and let $\xi= a^\infty\cdot b^\infty$. Clearly $\tra^{-n}\xi$ converges either to $a^\infty\cdot a^\infty$ or to $ b^\infty\cdot b^\infty$ if either $n\to \pm \infty$, in the topology of $\csp$. However, the orbit closure $\Xi_\xi=\overline{\Orb(\xi)}$ cannot be approximated in the Hausdorff topology by any periodic subshift. For indeed, $\dic(\xi)\cap \as^2=\{aa,ab,bb\}$. On the other hand, if $\eta$ is periodic of period $q$, there is a word of length $q$, say $u\in\as^q$, such that $\Orb(\eta)=\Orb(u^\infty)$. If $q=1$ then either $u=a^\infty$ or $u=b^\infty$. In either case $\dic(\eta)\cap \as^2$ does not contains the word $ab$. If $q>1$ then both letter $a,b$ occur in $u$ at least once, so that, by periodicity, $\dic(\eta)\cap \as^2$ contains the word $ba$. In both cases $\dic(\eta)\cap \as^2\neq \{aa,ab,bb\}$. Using the definition of dictionaries, it follows that $\dic(\eta)\cap \as^n\neq \dic(\xi)\cap \as^n$ for $n>1$. Hence there cannot be any periodic approximation of the subshift associated with $\xi$.
}%%
\hfill $\Box$
\end{exam}
%%%%%%%%%%

\noindent It is natural to ask if the condition of a subshift being topological transitive is closed in the Hausdorff topology on $\Inv$. The following example shows that this is not the case. Specifically, a sequence of periodic subshifts (being topological transitive) is defined and it is shown that its limit is not topological transitive.

%%%%%%%%%%%%%%%%
\begin{exam}
\label{1D.exam-TopTrans}
{\em Let $\as:=\{a,b\}$ and define the sequence of periodic elements $\eta_n\in\csp$ by $\eta_n:=(a^n b a^n bb)^\infty$. Clearly, their dictionary $\dic(\eta_n)\cap\as^m$ equals

$$
\{ a^m,ba^{m-1},aba^{m-2},\ldots,a^{m-1}b,bba^{m-2},abba^{m-3},\ldots,a^{m-2}bb\}
$$

\noindent for all $n>m$. Hence, the sequence of dictionaries $(\dic(\eta_n))_n\in\Dic$ converges in the local pattern topology to $\dic\in\Dic$ where $\dic$ is determined by $\dic\cap\as^m:=\dic(\eta_{m+1})\cap\as^m$ by Corollary~\ref{1D.cor-ConSubDic}. Denote by $\Xi\in\Inv$ the associated subshift with $\dic$. Due to Theorem~\ref{1D.th-DictHomeoSubs}, the subshifts $\Xi_n:=\Orb(\eta_n),\, n\in\NM,$ converge in the Hausdorff topology to $\Xi$. Assume $\Xi$ would be topological transitive, i.e., there is a $\xi\in\Xi$ such that $\overline{\Orb(\xi)}=\Xi$ implying $\dic(\xi)=\dic$ by Proposition~\ref{1D.prop-DicShi}. Then $b$ and $bb$ appear somewhere in $\xi$ as $b,bb\in\dic$. Thus, there is an $N\in\NM$ such that either $ba^Nbb$ or $bba^Nb$ appear in $\xi$. This is a contradiction as none of the words $ba^Nbb$ or $bba^Nb$ can be elements of $\dic$ by construction. Hence, $\Xi$ is not topological transitive.
}%%
\hfill $\Box$
\end{exam}
%%%%%%%%%%%%%%%%

%%%%%%%%%%%%%%%%%%%%%%%%%%%%%%%%%%%%%%%%%%%%%%%%%%%%%%%%%%%%%%%%%%%%%%%%%

\section{GAP-graphs and periodic approximations}
\label{1D.sect-GAP}

%%%%%%%%%%%%%%%%%%%%%%%%%%%%%%%%%%%%%%%%%%%%%%%%%%%%%%%%%%%%%%%%%%%%%%%%%

%%%%%%%%%%%%%%%%%%%%%%%%%%%%%%%%%%%%%%%%%%%%%%%%%%%%%%%%%%%%%%%%%%%%%%%%%

\subsection{Generalities on graphs}
\label{1D.ssect-GraphRem}

\noindent An {\em (oriented) graph} is a triple $\gs=(\vs,\es,\partial)$ where $\vs$ and $\es$ are discrete countable sets and $\partial_i$ (for $i=0,1$) are maps $\partial_i:\es\to\vs$ called {\em boundaries}. An element of $\vs$ is called a {\em vertex}, while an element of $\es$ is called an {\em edge}. An edge $e\in\es$ can be seen as an {\em arrow} joining its {\em origin} $\partial_0\,e$ to its {\em end} $\partial_1\,e$. An edge $e$ is called {\em outgoing} from the vertex $v$ if $\partial_0\,e=v$, while if $\partial_1\,e=v$ the edge $e$ is called {\em incoming}. The {\em vertex degree} ${\rm deg}(v)$ is defined by the number of incoming and outgoing edges. A vertex will be called {\em dandling} if it has either only incoming or only outgoing edges. If $v\in\vs$ is not dandling and ${\rm deg}(v)>2$ it will be called {\em branching}.

\vspace{.3cm}

\noindent A graph will be called {\em finite} if the vertex set $\vs$ and the edge set $\es$ are finite. A graph is called {\em simple} whenever for all edges $e$ the source $\partial_0 e$ and the range $\partial_1 e$ are not equal and if for any two vertices $u,v\in\vs$ there is at most one  edge with origin $u$ and end $v$. It will be called {\em semi-simple} if edges linking one vertex to itself are not excluded.

\vspace{.3cm}

\noindent A {\em path} is a finite sequence $\gamma=(e_1,e_2,\ldots, e_k)$ of edges such that $\partial_1\,e_i=\partial_0\,e_{i+1}$ for $1\leq i\leq k-1$. The number of edges $|\gamma|:=k$ in a path $\gamma$ is called the {\em length of $\gamma$}. The {\em origin} of a path $\gamma=(e_1,e_2,\ldots, e_k)$ is $\partial_0 \gamma:=\partial_0 e_1=u$, while its {\em end} is $\partial_1 \gamma:=\partial_1 e_k=v$. It will also be convenient to write a path as a map $\gamma:u\to v$. A path $\gamma=(e_1,\ldots,e_k)$ is said to visit the vertices $\{u_i=\partial_0 e_{i+1}\;:\; 1\leq i\leq k-1\}\subseteq\vs$ and $u_k=v=\partial_1\gamma$. Furthermore, a path $\gamma$ is called {\em closed} whenever $\partial_0 \gamma=\partial_1 \gamma$. A graph $\gs$ is called {\em connected}, if for any pair of distinct vertices $u,v\in\vs$ there is at least one path connecting them, namely $\gamma:u\to v$ or $\gamma:v\to u$.

%%%%%%%%%%%%%%%%
\begin{defini}
\label{1D.def-StrongConnect}
A graph $\gs$ is called {\em strongly connected}, if for any given pair $u,v\in\vs$ of vertices there are paths $\gamma,\gamma'$ such that $\gamma:u\to v$ and $\gamma':v\to u$. 
\end{defini}
%%%%%%%%%%%%%%%%

\noindent By definition a strongly connected graph is connected whereas the converse is in general false. Example~\ref{1D.exam-NotStrConn} provides an example of a connected but not strongly connected graph. A closed path $\gamma=(e_1,\ldots,e_k)$ is said to visit the vertices $\{u_i=\partial_0 e_i\;:\; 1\leq i\leq k\}\subseteq\vs$.

%%%%%%%%%%%%%%%%
\begin{proposi}
\label{1D.prop-ClosPathStrCon}
Let $\gs=(\vs,\es,\partial)$ be a finite graph. Then, there exists a closed path $\gamma$ that visit all vertices, if and only if $\gs$ is strongly connected.
\end{proposi}
%%%%%%%%%%%%%%%%

\noindent {\bf Proof: } This is clear from the definition.
\hfill$\Box$

\vspace{.3cm}

\noindent A closed path with the properties described in Proposition~\ref{1D.prop-ClosPathStrCon}, i.e. it visit all vertices, is called a \emph{global path} of $\gs$.

%%%%%%%%%%%%%%%%
\begin{defini}
\label{1D.def-Subgraph}
Let $\gs=(\vs,\es,\partial)$ be a graph. A subgraph $\gs'=(\vs',\es',\partial)$, denoted by $\gs'\preceq\gs$, is a graph such that (i) $\vs'\subseteq \vs$, (ii) $\es'\subseteq \es$ and (iii) every edge of $\gs'$ has its two boundaries in $\vs'$.
\end{defini}
%%%%%%%%%%%%%%%%

%%%%%%%%%%%%%%%%
\begin{defini}
\label{1D.def-gmap}
Given two graphs $\gs=(\vs,\es,\partial)$  and $\gs'=(\vs',\es',\partial')$, a graph map $\phi:\gs\to\gs'$ is a pair of maps $\phi=(\phi^{v}, \phi^{e})$, such that $\phi^{v}:\vs\to\vs'$ and $\phi^{e}:\es\to \es'$ and that $\partial' \phi^{e}(e)= \phi^{v}(\partial e)$ for all edge $e\in\es$. 

\noindent Given two graph maps $\gs\stackrel{\phi}{\to}\gs'\stackrel{\phi'}{\to}\gs"$, their composition $\phi'\circ\phi:\gs\to\gs"$ is defined by $\phi'\circ\phi=((\phi')^v\circ\phi^v,(\phi')^e\circ\phi^e)$. 
\end{defini}
%%%%%%%%%%%%%%%%

\noindent Clearly, the composition $\phi'\circ\phi$ is also a graph map. In addition, the identity map $id:\gs\to\gs$ is defined in the obvious way as well as the inverse of a graph map, if it is defined. A combinatorial graph is defined by $\vs=\{1,\cdots,N\}\subset \NM$, $\es\subseteq \vs^2\times \{0,\cdots,M\})$, where $\{0,\cdots,M\}\subset \NM_\ast =\NM\cup\{0\}$, with the convention that an edge $e= (i,j,0)$ does not exist and, for $m>0$, $(i,j,m)$ represent the $m$-th edge linking $i$ to $j$. Then the boundary operations are defined by $\partial_0 (i,j,m)=i\,,\, \partial_1 (i,j,m)=j$. Modulo graph isomorphisms, any finite graph is therefore equivalent to a combinatorial graph. In particular the set of finite graphs with $N$ vertices and $M$ edges is finite. A combinatorial graph is semi-simple if and only if $m\in\{0,1\}$ for any edge. It is simple if and only if, in addition, $(i,i,m)$ can only have $m=0$.

%%%%%%%%%%%%%%%%%%%%%%%%%%%%%%%%%%%%%%%%%%%%%%%%%%%%%%%%%%%%%%%%%%%%
 \subsection{GAP-graphs}
 \label{1D.ssect-GAP}

\noindent In 1894 Flye  introduced a graph representing possible continuation of words of finite lengths in \cite{Fl94}. In  1946 De Bruijn \cite{Bru46} and Good \cite{Goo46} specified independently the construction of these graphs. However, the name \emph{de Bruijn graphs} became common in a large community. These graphs encode the local structure of the one-dimensional Delone system according to the discussion given in Section \ref{1D.ssect-DelFLC}. Rauzy \cite{Rau83} provided the first use of the de Bruijn graphs in 1983 in order to compute the subword complexity. So these graphs are also called {\em Rauzy graphs} elsewhere \cite{Ca97,Ju10}. These graphs turn out to be nothing but the one-dimensional version of the Anderson-Putnam complex \cite{AP98} extended by G\"ahler \cite{Ga02} (never published). The notion proposed by G\"ahler can be found in \cite{Sa03,Sa08}. We expect that this complexes are crucial for the higher dimensional situation. Since our project aims at dealing eventually with higher dimensional systems, the Anderson-Putnam complex is the correct object to study. In view of the contribution of G\"ahler to the field, these graphs, defined below, will be called {\em GAP-graphs}. It is important to remark, though, that there are the same as the de Bruijn or Rauzy gaphs.

\vspace{.1cm}

\noindent In \cite{Mo05} the existence of so called de Bruijn sequences is studied for a given set of words of the same length. The concept of a strongly connected de Bruijn graph is already emphasized there. The equivalence between strongly connected de Bruijn graphs and irreducibility of the associated subshift of finite type is the most immediate consequence \cite[Lemma~9]{Mo05}. It is important to notice that the author of \cite{Mo05} uses a different notion of a dictionary than in the present work, though: a dictionary is only a set of words of a fixed length $n\in\NM$. However this result applies directly to the present context. The property of strongly connectedness is used in the work to construct periodic elements (Definition~\ref{1D.def-AssPerWor}). This idea of using de Bruijn graphs is not new \cite{Fi14}. However, the main motivation for the use of them here, came from the study of associated Schr\"odinger-like operators. The deep connection between these graphs and the Hausdorff topology on $\Inv$ ought to be underlined (see Theorem~\ref{1D.th-DictHomeoSubs}).

%%%%%%%%%%%%%%%%
\begin{defini}
\label{1D.def-GAP}
Let $\dic\in\Dic$ be a dictionary. For $k\in\NM$ define the vertex set $\vs_k:=\dic\cap\as^k$ and the edge set $\es_k:=\dic\cap\as^{k+1}$. The boundary maps $\partial_0\,,\,\partial_1:\es_k\to \vs_k$ are defined by

$$\partial_0(a_0a_1\ldots a_k):=
   a_0a_1\ldots a_{k-1}\,,
    \hspace{2cm}
     \partial_1(a_0a_1\ldots a_k):=
      a_1a_2\ldots a_k\,.
$$

\noindent The corresponding oriented graph $\gs_k:=\gs_k(\dic):=(\vs_k,\es_k,\partial)$ is called the {\em GAP-graph of $\dic$ of order $k$}. The sequence $\GAP:=\GAP(\dic):=(\gs_k)_{k\in\NM}$ will be called {\em GAP-sequence of the dictionary $\dic$}.
\end{defini}
%%%%%%%%%%%%%%%%

%%%%%%%%%%%%%%%%
\begin{rem}
\label{1D.rem-GAP}
{\em The boundary maps $\partial_0,\partial_1$ provide a constraint on the neighboring vertices and edges. As a matter of fact for edges $e=a_0a_1\ldots a_k,\, \te=b_0b_1\ldots b_k\in\es_k$ the condition $\partial_1 e=\partial_0 \te$ implies that $a_{i+1}=b_i$ for all $i=0,\ldots, k-1$. Thus, the GAP-graph encodes through the boundary maps what kind of continuation a word of length $k+1$ might have to the right and to the left. Specifically, if for some edge $e\in\es_k$ there is exactly one edge $\te$ fulfilling $\partial_0 e=\partial_1 \te$ (resp. $\partial_1 e=\partial_0 \te$) then the word of length $k+1$ associated with $e$ has a unique continuation to the left (resp. to the right). This observation emphasizes the importance of the branching vertices. For instance, in the case of  an outgoing branching vertex the word associated with the incoming edge does not have a unique continuation to the right.}%%
\hfill $\Box$
\end{rem}
%%%%%%%%%%%%%%%%
  
\noindent Since every subshift $\Xi\in\Inv$ is uniquely associated with a dictionary $\dic(\Xi)$ by Theorem~\ref{1D.th-DictHomeoSubs}, there is no ambiguity to use the notation $\gs_k(\Xi)$ for the GAP-graph of order $k\in\NM$ associated with $\dic(\Xi)$.

%%%%%%%%%%%%%%%%
\begin{proposi}
\label{1D.prop-GAPBasProper}
Let $\dic\in\Dic$ be the dictionary. For each $k\in\NM$ the GAP-graph $\gs_k:=\gs_k(\dic)$ is semi-simple with no dandling vertex. Furthermore, if $\dic=\dic(\xi)$ for some $\xi\in\as^{\ZM}$, then every GAP-graph $\gs_k$ of order $k\in\NM$ is connected.
\end{proposi}
%%%%%%%%%%%%%%%%

\noindent {\bf Proof: } Let $\dic\in\Dic$ be a dictionary with associated GAP-graph $\gs_k$ of order $k\in\NM$. Let $u\neq v$ be vertices in $\vs_k$ with an edge $e=a_0\ldots a_{k}\in\es_k$ connecting them, i.e. $u=a_0\ldots a_{k-1}=\partial_0e$ and $v=a_1\ldots a_{k}=\partial_1e$. In particular, any other edge linking $u$ to $v$ must be equal to $e$ meaning that $\gs_k$ is semi-simple. For every vertex $v=a_1\ldots a_k\in\vs_k$, there are letters $a,b\in\as$ such that $avb\in\dic$ by Definition~\ref{1D.def-dict}~(D.3). Thus, $v$ is not a dandling vertex.

\vspace{.1cm}

\noindent Let $\xi\in\csp$ be such that $\dic=\dic(\xi)$. Let $u,v$ be two distinct vertices $u\neq v\in\vs_k$. Then there is at least one word $w\in\dic(\xi)$ of length $|w|>k$ with subwords $u$ and $v$. In the language of the GAP-graphs this means the existence of at least one path, linking $u$ and $v$. Indeed, $u,v\in\vs_k$ occur in $\xi$, namely there is a $i,j\in\ZM$ such that $\xi_{[i,i+k-1]}=u$ and $\xi_{[j,j+k-1]}=v$. Without loss of generality, suppose $i<j$ and define the edges

$$e_l
	:=\xi_{[i+l,i+l+k]}
		\,,\qquad 0\leq l\leq j-i-1\,.
$$

\noindent By construction $\gamma=(e_0,\ldots,e_{j-i-1})$ is a path in $\gs_k$ satisfying $\partial_0 e_0=u$ and $\partial_1 e_{j-i-1}=v$. Hence, $\gs_k$ is connected.
\hfill $\Box$

\vspace{.2cm}

\noindent Note that $\dic=\dic(\xi)$ is equivalent to the fact that the associated subshift $\Xi:=\phi^{-1}(\dic)$ is topological transitive, i.e., $\Xi=\overline{\Orb(\xi)}$. Here, $\phi:\Inv\to\Dic$ denotes the homeomorphism defined in Theorem~\ref{1D.th-DictHomeoSubs}. The property that the GAP-graphs are strongly connected is crucial for periodic approximations by Theorem~\ref{1D.th-ExPerAppr} proven below. Before, a hereditary property is proven.

%%%%%%%%%%%%%%%%
\begin{proposi}
\label{1D.prop-HereGAP}
Let $\dic\in \Dic$ be a dictionary with GAP-sequence $\GAP(\dic):=(\gs_k)_{k\in\NM}$. If the GAP-graph $\gs_{k_0}$ of order $k_0\in\NM$ is strongly connected then every GAP-graph $\gs_k$ of order $k\leq k_0$ is strongly connected. In particular, if the GAP-graph $\gs_{l_0}$ of order $l_0\in\NM$ is not strongly connected, then all GAP-graphs $\gs_l$ of order $l\geq l_0$ are not strongly connected.
\end{proposi}
%%%%%%%%%%%%%%%%

\noindent {\bf Proof: } Let $\gs_{k_0}$ be the GAP-graph of order $k_0\in\NM$ in a GAP-sequence which is strongly connected. Let $k<k_0$ and $u,v\in\vs_k=\dic\cap\as^k$. As $\dic$ is a dictionary there are two words $\tilde{u}=a_1\ldots a_{k_0}$ and $\tilde{v}=b_1\ldots b_{k_0}$ contained in $\dic$ such that $a_1\ldots a_k=u$ and $b_{k_0-k+1}\ldots b_{k_0}=v$. Since $\tilde{u},\tilde{v}\in\vs_{k_0}$  there exists by hypothesis two paths $(\te_1,\ldots, \te_l)$ and $(e'_1,\ldots, e'_m)$ in $\gs_{k_0}$ such that $\tilde{u}=\partial_0 \te_1=\partial_1 e'_m$ and $\tilde{v}=\partial_1 \te_l=\partial_0 e'_1$. By convention, the representation $e=e(1)\ldots e(k_0+1)\in\gs_{k_0}$ is used. This paths induces two paths joining $u$ to $v$ respectively $v$ to $u$. For indeed, define $e_j:=\tilde{e}_1(j)\ldots\tilde{e}_1(j+k)$ for $j=1,\ldots, k_0-k+1$ and $e_{k_0-k+j}:=\tilde{e}_{j}(k_0-k+1)\ldots\tilde{e}_{j}(k_0+1)$ for $j=2,\ldots,l$. Then, $(e_1,\ldots, e_{k_0-k+l})$ defines a path in $\gs_k$ such that $\partial_0 e_1=u$ and $\partial_1 e_{k_0-k+l}=v$. Analogously, a path can be constructed joining $v$ to $u$ using the path $(e'_1,\ldots, e'_m)$ in $\gs_{k_0}$. Since $u,v\in\vs_k$ were arbitrary, the graph $\gs_k$ is strongly connected.
\hfill $\Box$

\vspace{.2cm}

\noindent The ``converse'' of Proposition~\ref{1D.prop-HereGAP} does not hold in general. Precisely, if the GAP-graph $\gs_{k_0}$ of order $k_0\in\NM$ is strongly connected the GAP-graphs of order $k>k_0$ might not be strongly connected as shown by the following example.

\begin{figure}[ht]
   \centering
\includegraphics[width=12cm]{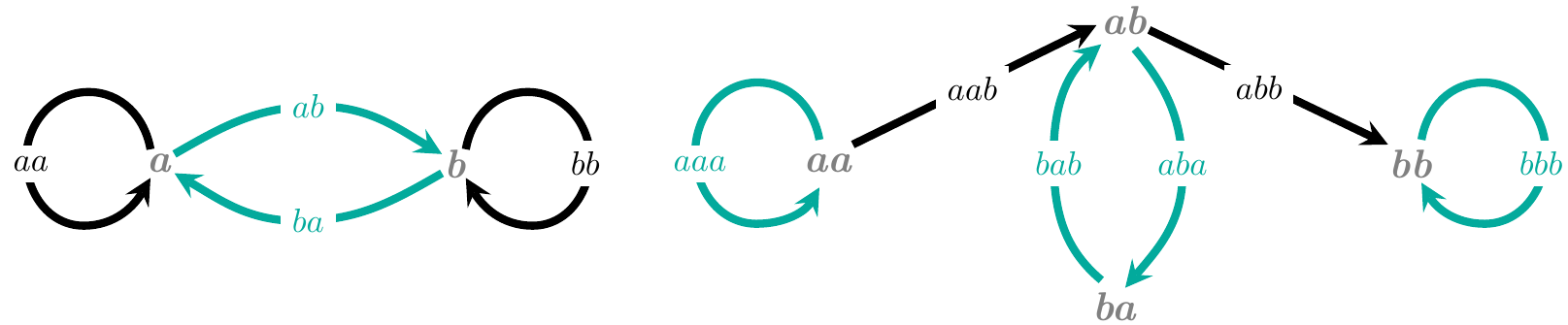}
\caption{The GAP-graph of order $1$ and $2$ of $\Xi$ defined in Example~\ref{1D.exam-NotStrConn}.}
\label{1D.fig-GAPNotStrongConn}
\end{figure}

%%%%%%%%%%%%%%%%
\begin{exam}
\label{1D.exam-NotStrConn}
{\em Let $\as:=\{a,b\}$ and let $\Xi:=\overline{\Orb(\xi)}$ be defines by $\xi:=a^\infty ba\cdot b^\infty\in\csp$. Then the GAP-graph $\gs_1$ of order $1$ is strongly connected while the GAP-graphs $\gs_k$ of order $k\geq 2$ are only connected but not strongly connected by Proposition~\ref{1D.prop-HereGAP}, c.f. Figure~\ref{1D.fig-GAPNotStrongConn}. Thus, $\Xi$ is not periodically approximable by Theorem~\ref{1D.th-ExPerAppr}.
}%%
\hfill $\Box$
\end{exam}
%%%%%%%%%%%%%%%% 

%%%%%%%%%%%%%%%%%%
%%%%%%%%%%%%%%%%%% 

\noindent The growth of the length of closed paths in a GAP-sequence is an indicator of the aperiodicity of the system, as shown in the following result. 

%%%%%%%%%%%%%%%%
\begin{proposi}
\label{1D.prop-ClosedPathBehaviour}
Let $\Xi\in\Inv$ and $\GAP=(\gs_k)_k$ be the corresponding GAP-sequence.
%%%%%%%%%%%%%%%%
\begin{itemize}
\item[(i)] If $\Xi$ is completely aperiodic then for each $m\in\NM$ there is a $k_0\in\NM$ such that for all $k\geq k_0$ the GAP-graph $\gs_k$ has no closed path of length smaller than or equal to $m$. In particular, the GAP-graphs are eventually simple.
\item[(ii)] Let $\Xi$ be not completely aperiodic, i.e. there is an $\eta\in\Xi$ with period $m$. Then, for each GAP-graph $\gs_k$ there exists a closed path of length $m$.
\item[(iii)] Let $\Xi=\Orb(\eta)$ be strictly periodic where $\eta$ has period $m\in\NM$. Then the GAP-graph $\gs_k$ has no branching vertex if $k\geq m-1$.
\end{itemize}
%%%%%%%%%%%%%%%%
\end{proposi}
%%%%%%%%%%%%%%%%

\noindent {\bf Proof: } Let $\Xi$ be a subshift, $m\in\NM$ and $\gs_k$ be the GAP-graph of order $k> m$. Let $\gs_k$ have a closed path $\gamma=(e_1,\ldots,e_m)$ with $e_1=a_0a_1\ldots a_{k}$. Due to the constraints given by the boundary maps  $\partial_1 e_i=\partial_0 e_{i+1}$ and $\partial_0 e_1=\partial_1 e_m$ it follows that this path corresponds to a long word $a_0a_1\ldots a_{k+m}$ with the condition $a_i = a_{i+m}$ for all $0\leq i< k$. These constraints imply that $e_1$ can be expressed as $e_1=v_ku$ for $v_k:=(a_0a_1\ldots a_{m-1})^{j_k}$ where $j_k\in\NM$ is the biggest number such that $j_k m\leq k$ and $u$ is a word of length smaller than $m$. By construction $u$ is a prefix of $v_k$. In addition $0\leq k- mj_k<m$ so that $k-m< |v_k|\leq k$. Since $e_1\in\es_k=\dic(\Xi)\cap\as^k$ it follows that  all the subwords of $e_1$ are contained in $\dic(\Xi)$. In particular $v_k\in\dic(\Xi)$, namely $\dic(\Xi)$ contains a word that is the $j_k$-th concatenation of a word of length $m$. 

\vspace{.1cm}

\noindent (i) Let $\Xi$ be completely aperiodic. By contradiction, let the GAP-sequence be such that, for each $k_0\in\NM$ there is a $k\geq k_0$ such that $\gs_k$ has a closed path of length $m$. By the previous considerations, there is a sequence of words $v_k\in\dic\,,\;k\in\NM\,,$ where each $v_k$ is the $j_k$-th time concatenation of a word of length $m$ and $\lim_{k\to\infty}|v_k|=\infty$. Since the alphabet is finite, there is a subsequence of words $v_{k_n}\in\dic$ which are the $j_{k_n}$-th time concatenation of a fixed word $u= a_0a_1\ldots a_{m-1}$ and $\lim_{n\to\infty}|v_{k_n}|=\infty$. By a standard argument, this implies that the periodic word $u^\infty$ is an element of $\Xi$, a contradiction as $\Xi$ is completely aperiodic. 

\vspace{.1cm}

\noindent Let $k_0\in\NM$ be chosen such that $\gs_k$ for $k\geq k_0$ does not contain a closed path of length $1$. Thus, the GAP graph $\gs_k$ is simple for every $k\geq k_0$ by Proposition~\ref{1D.prop-GAPBasProper}.

\vspace{.1cm}

\noindent (ii) Let $\eta\in\Xi$ be periodic with period $m$. Then, the dictionary $\dic(\eta)$ of $\eta$ is a subset of $\dic(\Xi)$. For $k\in\NM$ set $e_j:=\eta(j)\eta(j+1)\ldots\eta(j+k)$ with $1\leq j\leq m$. The collection of these edges defines a path in $\gs_k$. Due to the $m$-periodicity of $\eta$  it follows that $\partial_0 e_1=\partial_1 e_m$, namely $(e_1,\ldots,e_m)$ is a closed path in $\gs_k$.

\vspace{.1cm}

\noindent (iii) If $k\geq m-1$ every word of length $k$ has a unique continuation to the left and to the right. Thus, $\gs_k$ has no branching vertices.
\hfill $\Box$

%%%%%%%%%%%%%%%%
\begin{rem}
\label{1D.rem-NoBranch}
{\em The bound in Proposition~\ref{1D.prop-ClosedPathBehaviour} (iii) cannot be improved but is not necessarily sharp for all examples. In fact, let $\as:=\{a,b\}$. Consider the periodic words $\eta:=(baa)^\infty$ and $\xi:=(baab)^\infty$ in $\csp$. Then, $\eta$ has period $3$ and the GAP-graph $\gs_1$ of order $1$ admits a branching vertex. On the other hand, $\xi$ has period $4$ and the GAP-graph $\gs_2$ of order $2$ does not admit any branching vertex.}%%
\hfill $\Box$
\end{rem} 
%%%%%%%%%%%%%%%%

%%%%%%%%%%%%%%%%%%%%%%%%%%%%%%%%%%%%%%%%%%%%%%%%%%%%%%%%%%%%%%%%%%%%%%%%%

\subsection{Pruning and deriving}
\label{1D.ssect-ConstrGAP}

\noindent In this subsection the notion of a derived graph is introduced as well as the concept of pruning. This provides an algorithm to construct the GAP-graphs and so the periodic approximations of a subshift.

%%%%%%%%%%%%%%%%
\begin{defini}[Pruning]
\label{1D.def-PrunGraph}
Let $\gs=(\vs,\es,\partial)$ be a graph. A subgraph $\gs'\preceq\gs$ with  and $\gs'=(\vs',\es',\partial)$ is called the {\em pruned graph with respect to $(\vs',\es')$} .
\end{defini}
%%%%%%%%%%%%%%%%

%%%%%%%%%%%%%%%%
\begin{proposi}
\label{1D.prop-PrunFullShift}
Let $\gs^{full}_k:=(\vs^{full}_k,\es^{full}_k,\partial)$ be the GAP-graph of order $k\in\NM$ associated with the full shift $\csp$. For each subshift $\Xi\in\Inv$ and $k\in\NM$, the GAP-graph $\gs_k(\Xi)$ is the pruned graph of $\gs^{full}_k$ with respect to $\vs':=\vs^{full}_k\cap\dic(\Xi)$ and $\es':=\es^{full}_k\cap\dic(\Xi)$.
\end{proposi}
%%%%%%%%%%%%%%%%

\noindent {\bf Proof: } This is clear by definition of the GAP-graphs.
\hfill$\Box$

%%%%%%%%%%%%%%%%
\begin{defini}[Deriving]
\label{1D.def-DerivGraph}
Let $\gs=(\vs,\es,\partial)$ be an oriented graph. The {\em derived graph} $\delta\gs =(\delta\vs, \delta\es, \partial)$ is defined by
%%%%%%%%%%%%%%%%
\begin{itemize}
\item[(DG.1)] $\delta\vs:=\es$;

\item[(DG.2)] $\delta\es :=\{(e,e')\in\es\times \es\,:\, \partial_1 e=\partial_0 e'\}$;

\item[(DG.3)] If $(e,e')\in\delta\es$, then, $\partial_0 (e,e') :=e$ and $\partial_1 (e,e') :=e'$.
\end{itemize}
%%%%%%%%%%%%%%%%
\end{defini}
%%%%%%%%%%%%%%%%

\noindent By iterating the derivation of graphs leads to the following statement.

%%%%%%%%%%%%%%%%
\begin{proposi}
\label{1D.prop-PathGraphs}
Let $\gs=(\vs,\es,\partial)$ be a graph. For $n\in\NM$ the derived graph $\delta^n\gs$ is constructed as follows. Let $\Pp_n$ be the set of paths of length $n$ in $\gs$. Then, $\delta^n\vs=\Pp_n$, $\delta^n\es=\Pp_{n+1}$ and if $\gamma=(e_0,e_1,\ldots,e_n)\in\Pp_{n+1}$, then,
%%%%%%%%%%%%%%%%
$$
\partial_0\gamma=(e_0,e_1,\ldots,e_{n-1})\,,
\qquad\quad
    \partial_1\gamma =(e_1,e_2,\ldots,e_n)\,.
$$
%%%%%%%%%%%%%%%%
\end{proposi}
%%%%%%%%%%%%%%%%

\noindent {\bf Proof: } This can be proven by induction. The claim holds for $n=1$ just by definition. If the claim holds for $n$, then, $\delta^{n+1}\vs=\delta^n(\delta\vs)=\delta^n\es=\Pp_{n+1}$. 
Moreover, an edge in $\delta^{n+1}\es$ is a pair $(\gamma,\gamma')\in\Pp_{n+1}\times \Pp_{n+1}$ such that $\partial_1\gamma=\partial_0\gamma'$. This means that   $\gamma=(e_0,e_1,\ldots,e_n)\in\Pp_{n+1}$ and $\gamma' = (e_1, \ldots, e_n, e_{n+1})\in\Pp_{n+1}$ with $\partial_1 e_n= \partial_0 e_{n+1}$. 
Therefore the pair $(\gamma,\gamma')$ can be identified with the path $(e_0, e_1, \ldots, e_n, e_{n+1})\in\Pp_{n+2}$ 
and $\gamma=\partial_0(\gamma,\gamma')$ is nothing but the first $n+1$ edges of this path, while $\gamma'=\partial_1 (\gamma,\gamma')$ is given by the last $n+1$ edges.
\hfill $\Box$

\vspace{.2cm}

\noindent In the case of a GAP graph $\gs_k(\dic)=:\gs_k$ of order $k$ associated with a dictionary $\dic\in\Dic$, an edge of the derived graph $\delta\gs_k$ is given by  a pair of words of the form $u=a_0a_1\ldots a_k$ and $v=b_0b_1\ldots b_k$ such that $\partial_1 v=\partial_0v$. This implies that $v=a_1a_2\ldots a_k b_k\in\es_n$ see Remark~\ref{1D.rem-GAP}. However the word obtained by combining them, namely $u\vee v :=a_0a_1\ldots a_k b_k$, is not necessarily a word of the dictionary $\dic$. This leads to the following result.

%%%%%%%%%%%%%%%%
\begin{proposi}[Deriving and pruning]
\label{1D.prop-DerivPrun}
Let $\dic\in\Dic$ be a  dictionary. Then the GAP-graph $\gs_{k+1}$ of order $k+1$ coincides with the subgraph of $\delta\gs_k$ obtained by eliminating the edges $(u,v)\in\delta\es_k$ leading to words $u\vee v$ that are not in the dictionary $\dic$.
\end{proposi}
%%%%%%%%%%%%%%%%

\noindent {\bf Proof: } Let $k\in\NM$. An edge of $\gs_k$ is a word of the form $e=a_0a_1\ldots a_k$ and it is also a vertex of $\delta\gs_k$ by definition. Moreover, $u=\partial_0 e=a_0a_1\ldots a_{k-1}$ while $v=\partial_1 e=a_1a_2\ldots a_k$. Hence, a pair of such words $(u,v)$ becomes an edge of $\delta\gs_k$ if and only if $v= a_1a_2\ldots a_k b\in\dic$ for some $b\in\as$. In general the word $u\vee v=a_0a_1\ldots a_n b$ might not be in the dictionary $\dic$. Let then $\hat{\delta}\gs_k$ be the subgraph obtained by eliminating these edges. Consider the  graph map $\jmath=(\jmath_\vs,\jmath_\es):\hat{\delta}\gs_k\to\gs_{k+1}$  defined (for each $k$) by $\jmath_\vs(u) := u$ and $\jmath_\es\left(u,v\right) := u\vee v$. All the edges of the derived graph that are not in the dictionary $\dic$ are pruned. Thus, $\jmath:\hat{\delta}\gs_k\to\gs_{k+1}$ is a well-defined graph isomorphism. 
\hfill $\Box$

%%%%%%%%%%%%%%%%%%%%%%%%%%%%%%%%%%%%%%%%%%%%%%%%%%%%%%%%%%%%%%%%%%%%%%%%%

\subsection{Periodically approximable subshifts}
\label{1D.ssect-PropGAP}

\noindent In this section, the main Theorem~\ref{1D.th-ExPerAppr} is proven and sufficient conditions are provided for a subshift being periodically approximable.

%%%%%%%%%%%%%%%%
\begin{defini}
\label{1D.def-AssPerWor}
Let $\dic\in\Dic$ be a dictionary over the alphabet $\as$. Consider a closed path $\gamma:=(e_1,\ldots, e_{l})$ in the GAP-graph $\gs_k$ for some $k\in\NM$. An edge $e_j$ is identified by its $k+1$-letters, i.e. $e_j=e_j(0)\ldots e_j(k)$. The {\em associated periodic word $\eta:=\eta(\gamma)\in\csp$ with the closed path $\gamma$} is defined by $\eta:= \big(e_1(0)e_2(0)\ldots e_l(0)\big)^\infty$.
\end{defini}
%%%%%%%%%%%%%%%%

\noindent By construction, $\eta$ is periodic with period bounded by $l$ the number of edges in the path $\gamma$.

%%%%%%%%%%%%%%%%
\begin{lemma}
\label{1D.lem-PerCont}
Let $\dic$ be a dictionary and $\gamma:=(e_1,\ldots, e_{l})$ be a closed path in the GAP-graph $\gs_k(\dic)$ of order $k\in\NM$. Then the associated periodic word $\eta:=\eta(\gamma)\in\csp$ satisfies

$$\dic(\eta)\cap\as^k 
   = \{\partial_0 e_j \;:\; 1\leq j\leq l\}
    \subseteq\vs_k\,.
$$
\end{lemma}
%%%%%%%%%%%%%%%%

\noindent {\bf Proof: } Since $\gamma$ is a closed path the chains $\gamma_j:=(e_j, \ldots, e_l,e_1,\ldots,e_{j-1}),\; 1\leq j\leq l\,,$ define also a closed path in $\gs_k$. Denote by $\eta_j:=\eta(\gamma_j)$ the associated periodic word with the closed path $\gamma_j$ (Definition~\ref{1D.def-AssPerWor}) where $\eta=\eta_1$. Due to the periodicity of $\eta$ and $\eta_j$ and their definition, the equation $\tra^j\eta_1=\eta_{j+1}$ holds for $0\leq j\leq l-1$. Thus, Proposition~\ref{1D.prop-DicShi} implies $\dic(\eta)=\dic(\eta_j)$. Consequently, the inclusion $\{\partial_0 e_j \;:\; 1\leq j\leq l\}\subseteq \dic(\eta)\cap\as^k$ follows.

\vspace{.1cm}

\noindent For the converse inclusion, let $v\in\dic(\eta)\cap\as^k$. By construction there is a $1\leq j\leq l$ such that $v=\eta_j|_{[0,k-1]}$. Then $e_j=\eta_j|_{[0,k]}$ leads to $v=\partial_0 e_j$.
\hfill$\Box$

\vspace{.2cm}

\noindent The previous assertion implies that closed path in the GAP-graph of order $k$ give rise to periodic systems that only admit words of length $k$ appearing in the initial system $\Xi\in\Inv$. The main difficulty in defining periodic approximations is to avoid forbidden patterns of the original system by taking periodic boundary conditions. Thus, closed paths in the GAP-graphs are the right notion to do so. With this at hand, Theorem~\ref{1D.th-ExPerAppr} can be proved:

\vspace{.1cm}

\noindent \textbf{Proof of Theorem~\ref{1D.th-ExPerAppr}:} Let $\Xi\in\Inv$ be periodically approximable and $(\Xi_k)_k\subseteq\Inv$ be a sequence of periodic subshifts tending to $\Xi$. Thanks to Corollary~\ref{1D.cor-ConSubDic}, the convergence of the subshifts can be expressed in terms of the convergence of the associated dictionaries. Thus, without loss of generality it can be assumed that $\dic(\Xi_k)\cap\as^{k+1}=\dic(\Xi)\cap\as^{k+1}$ for each $k\in\NM$. Consequently, $\gs_k(\Xi_k)=\gs_k(\Xi)$. Let $k\in\NM$ and $\eta_k\in\Xi_k$ with period $l_k$. Hence, $\dic(\eta_k)=\dic(\Xi_k)$ holds by Proposition~\ref{1D.prop-DicShi} as $\Xi_k=\Orb(\eta_k)$. Then the chain $\gamma:=(e_1,\ldots,e_{l_k})$ defined by 

$$e_j
	:=\eta_k(j)\ldots\eta_k(j+k)
	\,,\qquad 1\leq j\leq l_k\,,
$$

\noindent is a closed path in the GAP-graph $\gs_k(\Xi_k)$ by the periodicity of $\eta_k$. Since, thanks to Lemma~\ref{1D.lem-PerCont},

$$\big\{\partial_0 e_j \;:\; 1\leq j\leq l_k\big\} =
	\dic(\eta_k)\cap\as^k=\dic(\Xi)\cap\as^k
$$

\noindent the path $\gamma$ is also a global path in $\gs_k(\Xi)$. Hence, $\gs_k(\Xi)$ is strongly connected by Proposition~\ref{1D.prop-ClosPathStrCon}. 

\vspace{.1cm}

\noindent Let $\Xi\in\Inv$ be such that there it has a subsequence of GAP-graphs that are all strongly connected. Proposition~\ref{1D.prop-HereGAP} yields that all GAP-graphs $\GAP=(\gs_k)_k$ are strongly connected. Thus, Proposition~\ref{1D.prop-ClosPathStrCon} assures the existence of a global (closed) path $\gamma_k:=(e_1,\ldots, e_{l_k})$ in $\gs_k$. Precisely, for each $v\in\vs_k$ there is a $1\leq j\leq l_k$ such that $\partial_0 e_j=v$. Let $\eta_k:=\eta(\gamma_k)$ be the associated periodic word of $\gamma_k$ defined in Definition~\ref{1D.def-AssPerWor}. Since $\gamma_k$ is a global path in $\gs_k$, the set $\dic(\eta_k)\cap\as^k$ equals to $\vs_k$ by Lemma~\ref{1D.lem-PerCont}. Since $\eta_k$ is periodic the associated periodic subshift $\Xi_k:=\Orb(\eta_k)$ satisfies 

$$\dic(\Xi_k)\cap\as^k
	=\dic(\eta_k)\cap\as^k
	= \vs_k 
	=\dic(\Xi)\cap\as^k\,.
$$

\noindent Thus, the periodic subshifts $(\Xi_k)_k$ converge to $\Xi$ in $\Inv$ by Corollary~\ref{1D.cor-ConSubDic}.
\hfill $\Box$

\vspace{.2cm}

\noindent It is essential to choose a global path in Theorem~\ref{1D.th-ExPerAppr} as can be seen by the following example.

%%%%%%%%%%%%%%%%%%
\begin{figure}[ht]
   \centering
\includegraphics[width=13cm]{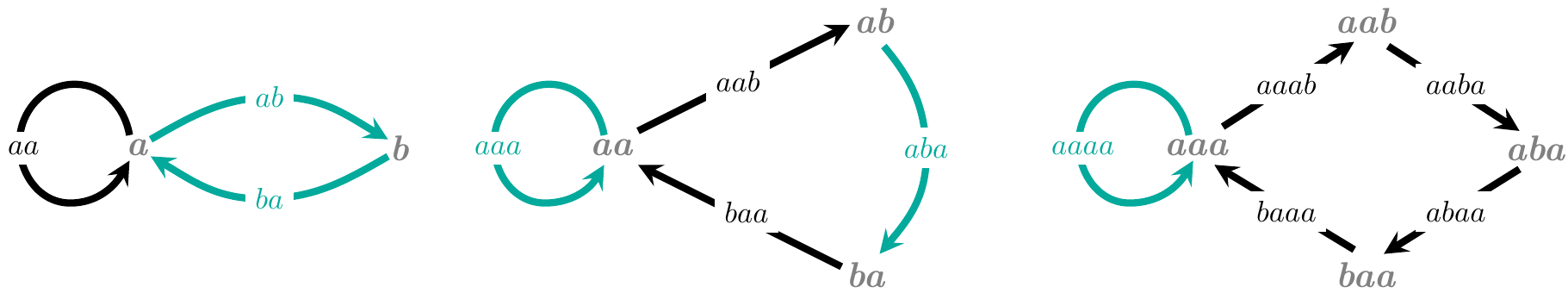}
\caption{The GAP-graph of order $1,2$ and $3$ for the One defect defined in Example~\ref{1D.exam-OneDef}.}
\label{1D.fig-GAPOneDef}
\end{figure}
%%%%%%%%%%%%%%%%%% 

%%%%%%%%%%%%%%%%
\begin{exam}
\label{1D.exam-OneDef}
{\em Let $\xi:= a^\infty b \cdot a^\infty$ be the one-defect sequence over the alphabet $\as:=\{a,b\}$. Then the associated subshift $\Xi:=\overline{\Orb(\xi)}$ has a family of strongly connected GAP-graphs by Theorem~\ref{1D.th-ExPerAppr} (see Figure~\ref{1D.fig-GAPOneDef}). The global paths are defined by the sequence $\eta_k:=(ba^k)^\infty$. The GAP-sequence $\GAP(\Xi)$ admits also a closed path $\gamma$ associated with the periodic word $\eta(\gamma)=a^\infty$ for every $\gs_k$ see Figure~\ref{1D.fig-GAPOneDef}. Clearly, the associated sequence of periodic subshifts $\Xi_k=\{\eta\}$ is constant and does not converge to $\Xi$. Hence, the choice of a global path is crucial in general. It also shows that $\Xi$ is not minimal.
}%%
\hfill $\Box$
\end{exam}
%%%%%%%%%%%%%%%% 

\noindent In contrast with Example~\ref{1D.exam-OneDef}, if $\Xi\in\Inv$ is minimal, every sequence of closed paths $\gamma_k$ in its GAP-graphs define a sequence of periodic subshifts converging to $\Xi$.

\vspace{.1cm}

\noindent {\bf Proof of Corollary~\ref{1D.cor-MinAllPath}:} Let $\Xi$ be minimal. According to \cite{RW92,La99A}, $\Xi$ is minimal if and only if each word occurs with bounded gaps, i.e for all $k\in\NM$ there exists a $l_k\in\NM$ such that any word $w\in\dic(\Xi)$ of length $|w|\geq l_k$ contains a copy of every elements of $\dic(\Xi)\cap\as^k$. Let $k\in\NM$, $u,v\in\dic(\Xi)\cap\as^k$ and $w=a_1\ldots a_{2l_k}\in\dic(\Xi)$. Thus, a copy of $u$ and $v$ appear in $w|_{[1,l_k]}$ and in $w|_{[l_k+1,2l_k]}$. Consequently, there are $w_1,w_2\in\dic(\Xi)$ such that $uw_1 v,vw_2 u\in\dic(\Xi)$ which are both subwords of $w$. In the same spirit as the proof of Proposition~\ref{1D.prop-HereGAP}, this two words give rise to two paths $\gamma_1$ and $\gamma_2$ in $\gs_k$ satisfying $\partial_0\gamma_1 = u = \partial_1 \gamma_2$ and $\partial_1\gamma_1 = v = \partial_0 \gamma_2$. Hence, $\gs_k$ is strongly connected.
\hfill$\Box$

\vspace{.2cm}

\noindent In the light of \cite{Mo05}, Lemma 9, $\Xi\in\Inv$ is periodically approximable if and only if all the associated subshifts of finite type are irreducible (namely, for each $u,v\in\dic(\Xi_k)$ there is a $w\in\dic(\Xi_k)$ such that $uwv\in\dic(\Xi_k)$). Specifically, the subshifts of finite type

$$\Xi_k
	:= \big\{ \xi\in\csp \;:\; \dic(\xi)\cap\as^k\subseteq\dic(\Xi)\cap\as^k \big\}
	\,,\qquad k\in\NM\,,
$$

\noindent are irreducible.

%%%%%%%%%%%%%%%%
\begin{coro}
\label{1D.cor-FinPerBruijn}
Let $\as$ be an alphabet. Then a subshift $\Xi\in\Inv$ can be approximated by a finite union of periodic subshifts if and only if there is a subsequence of GAP-graphs $(\gs_{k_l})_{l\in\NM}$ such that $\gs_{k_l}$ is a finite union of strongly connected graphs for every $l\in\NM$. 
\end{coro}
%%%%%%%%%%%%%%%%

\noindent Let $\Xi$ be a subshift induced by a substitution rule. Then, the substitution rule provides a method to build a subsequence of the GAP-graphs \cite{AP98}. Thus, a substitution rule gives an algorithm to compute periodic approximations. Recall that a substitution is a homomorphism $S:\as^\ast\to\as^\ast$ where $\as^\ast:=\dic(\csp)$. Furthermore, a substitution is called {\em primitive} if there is an $N\in\NM$ such that for any two letters $a,b\in\as$, the letter $a$ appears in $S^N(b)$. A primitive substitution defines uniquely a subshift $\Xi\in\Inv$ \cite{Qu87,Fog02,Qu10}.

\vspace{.1cm}

\noindent {\bf Proof of Proposition~\ref{1D.prop-PrimSub}: } Let $S$ be a primitive substitution with associated subshift $\Xi$. Since $\Xi$ is minimal \cite{Fog02}, $\Xi$ is periodically approximable by Corollary~\ref{1D.cor-MinAllPath} and $\gs_k(\Xi)$ is strongly connected for each $k\in\NM$ by Theorem~\ref{1D.th-ExPerAppr}. Consider a closed path $\gamma$ in the GAP-graph $\gs_l(\Xi)$ of order $l\in\NM$. Let $\eta:=\eta(\gamma)$ be the associated periodic word, c.f. Definition~\ref{1D.def-AssPerWor}. Lemma~\ref{1D.lem-PerCont} implies $\dic(\eta)\cap\as^l\subseteq\dic(\Xi)\cap\as^l$. Since $S$ is primitive, $\lim_{k\to\infty}|S^n(a)|=\infty$ follows for every letter $a\in\as$. Thus by the previous two considerations, there is an $k_0(m)\in\NM$ for every $m\in\NM$ such that 

$$\dic(S^k(\eta))\cap\as^m
	\subseteq\dic(\Xi)\cap\as^m
	\,,\qquad
	k\geq k_0(m)\,.
$$

\noindent Furthermore, for every $v\in\dic(\Xi)$, there exists a $k_v\in\NM$ such that $v$ is a subword of $S^{k_v}(a)$ for each $a\in\as$ as $S$ is primitive \cite{Fog02}. Consequently, for each $m\in\NM$ there is an $k_1(m)\in\NM$ with $k_1(m)\geq k_0(m)$ satisfying 

$$\dic(S^k(\eta))\cap\as^m
	=\dic(\Xi)\cap\as^m
	\,,\qquad
	k\geq k_1(m)\,.
$$

\noindent Hence, Corollary~\ref{1D.cor-ConSubDic} implies $\Xi_k\to\Xi$ where $\Xi_k=\Orb(S^k(\eta))$ as $\dic(S^k(\eta))=\dic(\Xi_k)$ by Proposition~\ref{1D.prop-DicShi}.
\hfill$\Box$

%%%%%%%%%%%%%%%%
\begin{rem}
\label{1D.rem-PrimSub}
{\em Proposition~\ref{1D.prop-PrimSub} extends to more general substitutions defining a subshift in the following way: Given a substitution over the alphabet $\as$ which is not necessarily primitive such that for at least one letter $a\in\as$, there is an $k_a\in\NM$ where $S^{k_a}(a)$ contains all letters in $\as$ and $uav=S^{k_a}(a)$ where $u,v$ are words that are non-empty. The latter conditions is needed so that $\dic(S)$ defined in the following satisfies Definition~\ref{1D.def-dict}~(D.3). Define $\dic(S)$ by all possible subwords appearing in any power $S^k(a)$ which defines a dictionary. Denote by $\Xi(S)\in\Inv$ the associated subshift of $\dic(S)$. If there is an $l\in\NM$ such that $\gs_l(\Xi(S))$ admits a closed path $\gamma$ such that the associated periodic word $\eta=\eta(\gamma)$ contains the letter $a$, then $S^k(\eta)$ defines a sequence of periodic subshifts $\Xi_k:=\Orb(S^k(\eta))$ converging to $\Xi$ in the Hausdorff topology on $\Inv$. The proof is similar to the proof of Proposition~\ref{1D.prop-PrimSub}. That no forbidden pattern is created is guaranteed by taking a closed path. That all possible words appear eventually follows by the fact that the letter $a$ appears in $\eta$.
}%%
\hfill $\Box$
\end{rem}
%%%%%%%%%%%%%%%%

\noindent Recall that the period of the elements of $\Xi_k$ is given by $q_k:=\sharp\Xi_k$. The period of the periodic approximations can be bounded from below by the quantity $\sharp\dic(\Xi)\cap\as^k$ which is called {\em subword complexity} of $\Xi$. This quantity is studied in more detail in Subsection~\ref{1D.ssect-defect}

%%%%%%%%%%%%%%%%
\begin{coro}
\label{1D.cor-PerGrowth}
Let $\Xi\in\Inv$ be periodically approximable and $(\Xi_k)$ be a sequence of periodic subshifts tending to $\Xi$. If $\Xi$ contains a non-periodic element, then the period $q_k=\sharp\Xi_k$ goes to infinity if $k\to\infty$.
\end{coro}
%%%%%%%%%%%%%%%%

\noindent {\bf Proof: } The period of $\eta\in\csp$ is given by the number $q$ of elements in the orbit $\Orb(\eta)$. In addition, because of the periodicity of $\eta$, the number of distinct words in $\dic(\eta)\cap\as^k$ cannot exceed the period of $\eta$. Hence $\sharp\dic(\eta)\cap \as^k\leq q$ for all $k\in\NM$. On the other hand, distinct element of $\dic(\eta)\cap\as^k$ give rise to a distinct element in $\Orb(\eta)$, leading also to $q=\sharp\Orb(\eta)\geq\sharp\dic(\eta)\cap\as^k$. Hence, if $\eta_k$ arises by a global path in $\gs_k=\gs_k(\Xi)$ then 

$$q_k
   := \sharp\Orb(\eta_k)
    \geq \sharp \dic(\eta)\cap\as^k
     = \sharp \dic(\Xi)\cap\as^k\,.
$$

\noindent According to \cite{HM40}, the limit $\lim_{k\to\infty}\sharp\dic(\xi)\cap\as^k$ goes to infinity, if and only if $\xi$ is non-periodic. Thus, $q_k$ tends to infinity if $\Xi$ contains a non-periodic element since $\dic(\xi)\cap\as^k\subseteq\dic(\Xi)\cap\as^k$ holds for $\xi\in\Xi$.
\hfill$\Box$

%%%%%%%%%%%%%%%%%%%%%%%%%%%%%%%%%%%%%%%%%%%%%%%%%%%%%%%%%%%%%%%%%%%%%%%%%

\subsection{Subword complexity function and defects}
\label{1D.ssect-defect}

%%%%%%%%%%%%%%%%%%%%%%%%%%%%%%%%%%%%%%%%%%%%%%%%%%%%%%%%%%%%%%%%%%%%%%%%%

\noindent This section introduces the concept of subword complexity function and investigates its connection with the existence of branching vertices in the GAP-graphs. As discussed in Remark~\ref{1D.rem-GAP} the branching vertices are related to the non-periodicity of a two-sided infinite word. It will be argued that the branching vertices play a role similar to boundary conditions in difference or differential equations. The concept of {\em amenability} of a subshift will follow.

\vspace{.1cm}

\noindent Let $\dic\in\Dic$ be a dictionary. The {\em subword complexity function} $\scf_\dic:\NM\to\NM$ is defined by $\scf_\dic(k):=\sharp\dic\cap\as^k$. If there is no ambiguity the notation $\scf:=\scf_\dic$ will be used instead. Due to Theorem~\ref{1D.th-DictHomeoSubs}, the subword complexity function is also defined for every subshift $\Xi\in\Inv$.

\vspace{.1cm}

\noindent The elementary estimate $\scf(k+1)\leq \sharp\as\,\scf(k)$ holds and $\scf(1)=\sharp\as$, leading to $\scf(k)\leq (\sharp\as)^k$. This upper bound corresponds exactly to the complexity function of the full shift $\as^\ZM$. Furthermore, the subword complexity function is bounded if the associated subshift $\Xi:=\phi^{-1}(\dic)$ is finite where $\phi:\Inv\to\Dic$ is the homeomorphism defined in Theorem~\ref{1D.th-DictHomeoSubs}. On the contrary, if there is a non-periodic element in $\Xi$, $\scf$ grows at least linearly, namely $\scf(k)\geq k+1$ \cite{HM40}. This function has been studied by many authors \cite{ER82,Rau83,AR91, Al94,Ca97,La99A,Ju10}. The following definition is inspired from \cite{AR91}

%%%%%%%%%%%%%%%%
\begin{defini}
\label{1D.def-BranchPoint}
Let $\Xi$ be a subshift with associated GAP-graphs $\GAP=(\gs_k)_k$. Let $u$ be a vertex of $\gs_k$. Then 

(i) $\partial^+ u$ ({\em resp.} $\partial^-u$) denotes the number of oriented edges $e$ starting from ({\em resp.} ending to) $u$ namely with $\partial_0 e=u$ ({\em resp.} $\partial_1 e=u$),

(ii) the number $\nbr^+(k)$ ({\em resp.} $\nbr^-(k)$) of forwards ({\em resp.} backward) branching vertices in $\gs_k$ is the number of vertices $u$ such that $\partial^+u>1$ ({\em resp.} $\partial^-u>1$),

(iii) the number $\nbr(k)$ of branching vertices in $\gs_k$ is defined as the number of $u$'s such that at least one of $\partial^\pm u$ is larger than $1$. 
\end{defini}
%%%%%%%%%%%%%%%% 

%%%%%%%%%%%%%%%%
\begin{proposi}
\label{1D.prop-branch}
Let $\as$ be a finite alphabet with more than one letter and let $\dic\in\Dic$. Then for any $k\in\NM$, $\nbr^\pm(k)\leq \scf(k+1)-\scf(k)$. In addition

$$\frac{\scf(k+1)-\scf(k)}{\sharp\as-1}
	\leq \nbr(k)
	\leq 2 (\scf(k+1)-\scf(k))
$$

\noindent The previous estimated are optimal.
\end{proposi}
%%%%%%%%%%%%%%%%

\noindent {\bf Proof: }The number of vertices of the GAP-graph $\gs_k$ is exactly the number of word of length $k$ in $\dic$, namely $\scf(k)$, while the number of edges is exactly the number of words of length $k+1$, namely $\scf(k+1)$. Let $L_v^+(j)$ denote the number of vertices with exactly $j$ outgoing edges. If $\as$ has exactly $A$ letters then $L_v^+(j)=0$ if $j>A$. Since there is no dandling vertex, each vertex admits at least one outgoing edge. Similarly, each edge has an origin, at some vertex. Consequently the total number of vertices and of edges are given by 

$$\sum_{j=1}^A L_v^+(j)=\scf(k)\,,
   \hspace{2cm}
    \sum_{j=1}^A j L_v^+(j)=\scf(k+1)\,.
$$

\noindent A vertex is forward branching if it has at least two outgoing edges. Hence 

$$\nbr^+(k)=\sum_{j=2}^A L_v^+(j)
   \leq \sum_{j=1}^A (j-1)L_v^+(j)=\scf(k+1)-\scf(k)\,.
$$

\noindent The same estimate can be obtained for the backward branching vertices leading to the upper bound for the total number of branching vertices as $\nbr(k)\leq \nbr^-(k) + \nbr^+(k)$. 

\vspace{.1cm}

\noindent For each forward branching vertex $u$, $\partial^+ u$ is bounded by $A$. Hence

$$\scf(k)-\nbr^+(k) + A\, \nbr^+(k)\geq \scf(k+1)
$$

\noindent follows implying the lower bound for the total number of branching vertices as $\nbr^+(k)\leq \nbr(k)$.

\vspace{.1cm}

\noindent The optimality of the estimates can be seen for an alphabet with two letters $a,b$: Consider the Fibonacci subshift $\Xi\in\Inv$ defined Subsection~\ref{1D.ssect-Fib}. For the Fibonacci substitution the subword complexity is given by $\scf(k)=k+1$ \cite{Fog02}. From the previous considerations, the estimates $1\leq \nbr(k)\leq 2$ hold. The GAP-graph $\gs_2$ of order $2$ of $\Xi$ has two branching vertices whereas the GAP-graph $\gs_3$ of order $3$ admits only one branching vertex see Figure~\ref{1D.fig-Fib}. Thus, the estimates on $\nbr:\NM\to\NM$ are optimal.
\hfill $\Box$

\vspace{.2cm}

\noindent It is known that if the subword complexity function satisfies $\scf(k)\leq k$ for $k$ large enough, the corresponding subshift is automatically eventually periodic periodic \cite{CH73,Fog02}.

%%%%%%%%%%%%%%%%
\begin{defini}
\label{1D.def-sturm}
A minimal subshift $\Xi\subset \as^\ZM$ is called {\em Sturmian}, whenever $\scf(k)=k+1$.
\end{defini}
%%%%%%%%%%%%%%%%

\noindent A Sturmian subshift can only be defined on an alphabet with two letters \cite{HM40}. In view of the comments made previously they are the aperiodic subshift with minimal subword complexity functions. Such subshifts are completely classified \cite{Fog02}. They correspond exactly to the 1D-tilings obtained by the cut-and-project method from $\ZM^2$ onto $\RM$.

%%%%%%%%%%%%%%%%
\begin{proposi}
\label{1D.prop-sturm}
The GAP-sequence of a Sturmian subshift is made of graphs with exactly one forward and one backward branching vertex, which may either be distinct or may coincide. 
\end{proposi}
%%%%%%%%%%%%%%%%

\noindent {\bf Proof: }Since $\scf(k+1)-\scf(k)=1$ this follows immediately from Proposition~\ref{1D.prop-branch}.
\hfill $\Box$

\vspace{.2cm}

\noindent It is also worth noticing that Arnoux and Rauzy \cite{AR91} have completely classified geometrically the subshifts with subword complexity functions $\scf(k)= (l-1)k+1$ for $l\geq 2$. On the other hand the complexity function for the Golay-Rudin-Shapiro sequence is $\scf(k)=8k-8$ for $k\geq 8$ \cite{AS93}.

%%%%%%%%%%%%%%%%
\begin{coro}
\label{1D.cor-ApBra}
If $\Xi\in\Inv$ is aperiodic, then all the GAP-graphs admit at least one branching vertex.
\end{coro}
%%%%%%%%%%%%%%%%

\noindent  {\bf Proof: } If $\Xi$ is aperiodic, then there is a non-periodic $\xi\in\Xi$. Thus, $\scf_{\dic(\xi)}$ grows strictly, namely $\scf_{\dic(\xi)}(k+1)>\scf_{\dic(\xi)}(k)$ \cite{HM40}. Hence, every GAP-graph of $\dic(\xi)$ admits at least one branching vertex by Proposition~\ref{1D.prop-branch}. Since these GAP-graphs are subgraphs of the GAP-graphs of $\Xi$ the desired result follows.
\hfill $\Box$

\vspace{.2cm}

\noindent Why are branching vertices relevant ? One possible way to understand this point is to consider the discrete Schr\"odinger equation as the eigenvalue equation for a self-adjoint operator $H$ of the form given in Equation~\ref{1D.eq-schdis}.

 (i) In some vague sense, the GAP-graphs $\gs_k$ could be seen as a finite volume approximation of the subshift.

 (ii)  If the word length $k$ is very large, much larger than the range of the operator $H$, the latter can be seen as a discrete operator acting on the Hilbert space $\ell^2(\vs_k)$, where $\vs_k$ is the set of vertices of the graph $\gs_k$. However, at each branching vertex, like in quantum graphs \cite{BK13}, a boundary condition has to be defined. In other words, branching vertices of the graph are acting as a boundary condition. 

 (iii) It can be guessed also, that if the boundary condition do not correspond to a closed cycle along the graph, it might create boundary states, looking like defects in the material \cite{Bel13}.

 (iv) Hence if the size of the boundary grows too fast, relative to the side of the system, namely here the number of vertices in $\gs_k$, as $k\to\infty$, the defects due to the branching vertices can become dominant in the spectrum of $H$. In Statistical Mechanics, such a growth occurs in models on an hyperbolic space, with non-negligible boundary effect at infinite volume. In Group Theory, this is related to the concept of amenability \cite{Gr69}.

\noindent This discussion leads to the following definitions

%%%%%%%%%%%%%%%%
\begin{defini}
\label{1D.def-amenability}
The subshift $\Xi$ will be called amenable whenever the subword complexity function satisfies 

$$\lim_{k\to\infty}
    \frac{\scf(k+1)}{\scf(k)}=1\,.
$$

\noindent The configurational entropy of $\Xi$ is defined as 

$$h= \limsup_{k\to\infty}
   \frac{\ln\{\scf(k)\}}{k}\,.
$$
\end{defini}
%%%%%%%%%%%%%%%%

\noindent The definition of amenability is leading to an asymptotic  negligible number of branching vertices when compared to the number of vertices in the GAP-graphs. On the other hand a nonzero configurational entropy gives asymptotically a finite proportion of branching vertices. The following result is immediate and its proof will be left to the reader.

%%%%%%%%%%%%%%%%
\begin{lemma}
\label{1D.lem-zeroh}
Let $\Xi$ be a subshift.

\noindent (i) If $\Xi$ is amenable its configurational entropy vanishes. 

\noindent (ii) If the configurational entropy $h$ is positive, then the proportion of branching vertices along the sequence of GAP-graphs is asymptotically bounded from below by $\left(e^{h}-1\right)/\sharp\as$.
\end{lemma}
%%%%%%%%%%%%%%%%

\noindent This discussion leads to two problems

%%%%%%%%%%%%%%%%
\begin{prob}[Structural Defects]
\label{1D.prob-defects}
{\em Is it correct to link the branching vertices in the sequence of GAP-graphs to the appearance of defects structurally created by the subshift itself~? In particular, would it be possible to establish, from a thermodynamical approach in Statistical Mechanics, that such defects indeed occur in a material due to the overall structure~? Such defects seem to be present in 3D-quasicrystalline alloys in particular as a contribution to the Density of State near the Fermi level, where a pseudo-gap, usually explained by a Hume-Rothery mechanism, is taking place, and partially filled (see a detailed discussion in \cite{Bel03}, Section 6.5). Physicists have long argued about whether such defects are coming from alien impurities introduced during the sample production or whether they occur spontaneously due to the internal structure of the material.
}%%
\hfill $\Box$
\end{prob}
%%%%%%%%%%%%%%%% 

%%%%%%%%%%%%%%%%
\begin{prob}[Nature of the Spectral Measure]
\label{1D.prob-spe}
{\em If the branching vertices are present in overwhelming numbers, can one expect the corresponding defect to

(a) fill the spectral gaps, at least near the energies influenced by the defects~?

(b) to create enough interferences in the wave functions to localize the quantum particles described by the Hamiltonian~? In other words, can a pure point spectrum result from too many such defects~? 
}%%
\hfill $\Box$
\end{prob}
%%%%%%%%%%%%%%%% 

\noindent An evidence in favor of the last question is provided by the Anderson model on $\ZM$, with a random potential satisfying a Bernoulli distribution, which has been proved to have a pure point spectrum \cite{CKM87}. Such a random potential can be described by a full shift, for which all vertices in the GAP-graphs are branching. 

\vspace{.3cm}

%%%%%%%%%%%%%%%%%%%%%%%%%%%%%%%%%%%%%%%%%%%%%%%%%%%%%%%%%%%%%%%%%%%%
\section{Two Examples}
\label{1D.sect-2ex}

\noindent To finish this article, two standard examples are described in more detail in this section. The first one is provided by the Fibonacci sequence which represents a paradigm for one-dimensional quasicrystals. The second the Golay-Rudin-Shapiro sequence which is still largely a mystery as far as the corresponding Schr\"odinger operators are concerned. It is believed that the distribution of letters is disordered enough to produce some point spectrum while it has the lowest possible algorithmic complexity and zero entropy. For this reason, it represents a borderline case separating continuous spectra from point spectra in one-dimension. The description of its GAP-graph given here in Fig.~\ref{1D.fig-GAPgrs} seems to be the first published one, as far as the authors have been able to check.

%%%%%%%%%%%%%%%%%%%%%%%%%%%%%%%%%%%%%%%%%%%%%%%%%%%%%%%%%%%%%%%%%%%%%%%%%
\subsection{The Fibonacci sequence}
\label{1D.ssect-Fib}

\noindent A typical example of an aperiodic subshift $\Xi$ is provided by the {\em Fibonacci sequence}. In this case the substitution $S$ is defined by $S:a\to ab\,,\,b\to a$. 

%%%%%%%%%%%%%%%%%%
\begin{figure}[ht]
   \centering
\includegraphics[width=13cm]{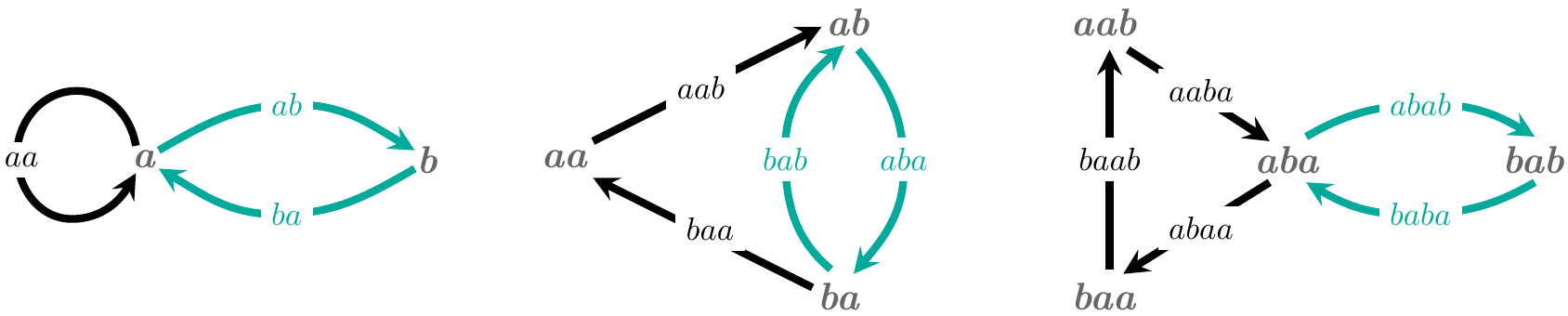}
\caption{The GAP-graph of order $1,2$ and $3$ for the Fibonacci subshift.}
\label{1D.fig-Fib}
\end{figure}
%%%%%%%%%%%%%%%%%% 

%%%%%%%%%%%%%%%%
\begin{proposi}
\label{1D.prop-Fib}
Let $\Xi$ be the Fibonacci subshift. For $k\in\NM$, let $\eta^a_k:=S^k(a^\infty)$ and $\eta^b_k:=S^k(b^\infty)$ be periodic sequences with corresponding periodic subshifts $\Xi^a_k:=\Orb(\eta^a_k)$ and $\Xi^b_k:=\Orb(\eta^b_k)$. Then $\Xi$ is periodically approximable and the sequences of periodic subshifts $(\Xi^a_k)_k$ and $(\Xi^b_k)_k$ converge to $\Xi$.
\end{proposi}
%%%%%%%%%%%%%%%%

\noindent {\bf Proof: } Define the closed path $\gamma=(e)$ in $\gs_1$ where $e$ is the edge $aa$, c.f. Figure~\ref{1D.fig-Fib}. Furthermore, $\eta^b_k=\eta^a_{k-1}$ holds for $k\geq 2$. Thus, the statement follows from Proposition~\ref{1D.prop-PrimSub} since the substitution $S$ is primitive.
\hfill $\Box$

\vspace{.2cm}

\noindent Combining Proposition~\ref{1D.prop-Fib} with Theorem~\ref{1D.th-ContSpectrGenSchrDyn} leads to

%%%%%%%%%%%%%%%%
\begin{coro}
\label{1D.cor-Fib}
Let $\Xi$ be the Fibonacci subshift and $H$ be a generalized discrete Schr\"odinger operator defined in Equation~\ref{1D.eq-schdis} satisfying (R1, R2, R3). Then the equations

$$\sigma(H_\xi)
	=\lim_{k\to\infty} \sigma(H_{\eta^a_k})
	=\lim_{k\to\infty} \sigma(H_{\eta^b_k})
$$

\noindent hold for every $\xi\in\Xi$.
\end{coro}
%%%%%%%%%%%%%%%%

\noindent This result proves that the numerical calculation by \cite{OK85} is justified to compute the spectrum of the Fibonacci Hamiltonian.

%%%%%%%%%%%%%%%%%%%%%%%%%%%%%%%%%%%%%%%%%%%%%%%%%%%%%%%%%%%%%%%%%%%%%%%%%
\subsection{The Golay-Rudin-Shapiro sequence}
\label{1D.ssect-GRS}

\noindent This sequence was defined and used in \cite{Sh51,Go51,Rud59} (see also \cite{AS93,Al97}). Let $n\in\NM$ be written in base $2$. $n=\epsilon_0+2\epsilon_1+\cdots+ 2^k\epsilon_k$ with $\epsilon_i\in\{0,1\}$. Then set $a_n=(-1)^{\epsilon_n\epsilon_{n+1}}\in\{+1,-1\}$. On the $2$-letters alphabet $\as=\{a,b\}=\{+1,-1\}$ it gives the sequence (the splitting is provided to make the reading easier)

$$\cdots aaab|aaba|aaab|bbab|aaab|aaba|bbba|aaba|aaab|aaba|aaab|bbab|bbba \cdots \,.
$$

\noindent It satisfies the following recursion formula
 
$$a_{2n} = a_n\,,
  \hspace{2cm}
   a_{2n+1} = (-1)^n a_n\,.
$$

\noindent A new alphabet $\bs=\{A,B,C,D\}$ is defined by setting $A=aa,B=ab,C=ba,D=bb$, so that this sequence is generated by the following substitution

$$S: A\to AB\,,
   \hspace{.7cm}
    B\to AC\,,
     \hspace{.7cm}
      C\to DB\,,
       \hspace{.7cm}
        D\to DC\,.
$$

%%%%%%%%%%%%%%%%%%
\begin{figure}[ht]
   \centering
\includegraphics[width=13cm]{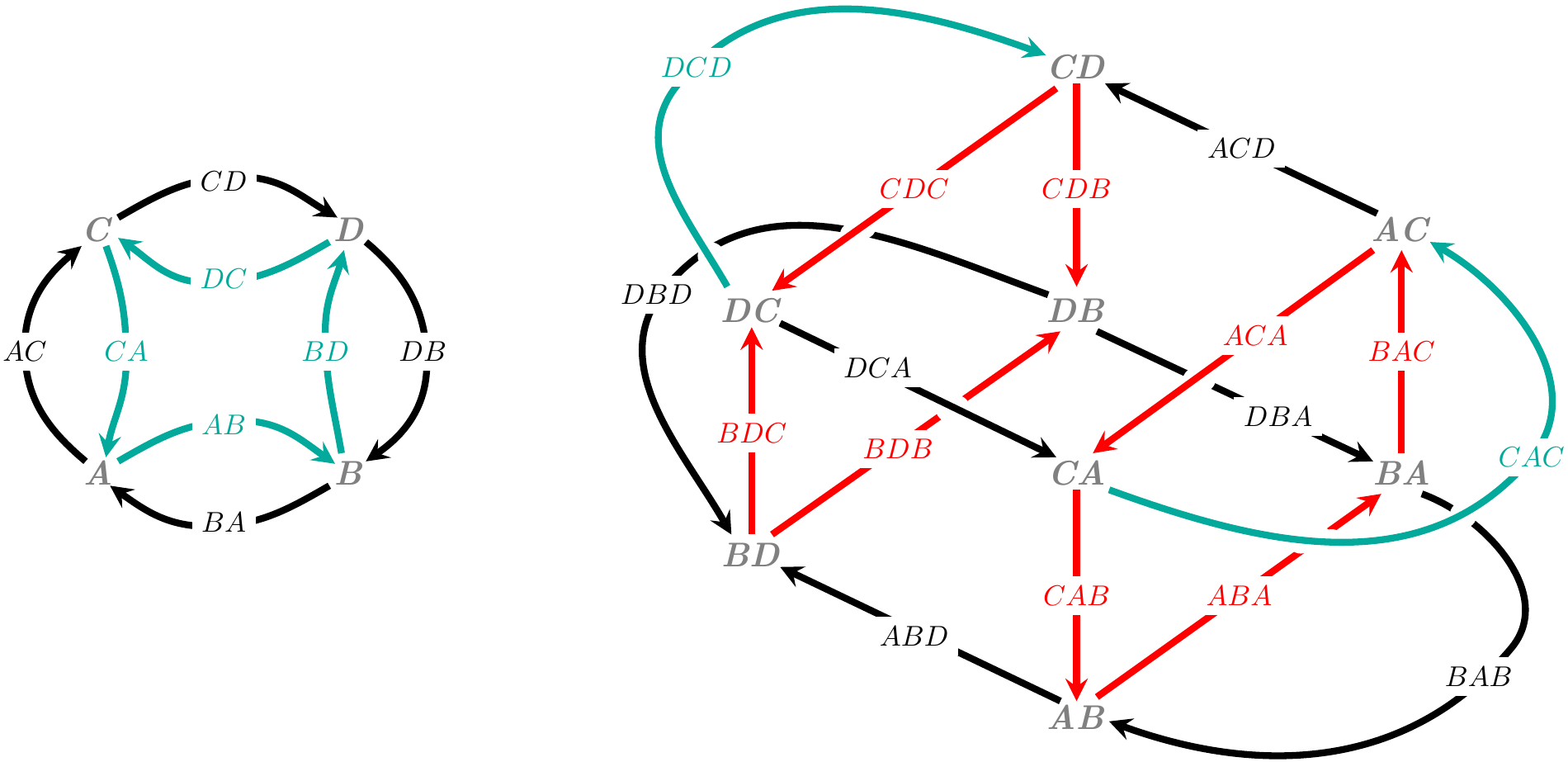}
\caption{The GAP-graph of order $1$ and $2$ for the Golay-Rudin-Shapiro subshift.}
\label{1D.fig-GAPgrs}
\end{figure}
%%%%%%%%%%%%%%%%%% 

%%%%%%%%%%%%%%%%
\begin{proposi}
\label{1D.prop-GRS}
Let $\Xi$ be the Golay-Rudin-Shapiro subshift and $x\in\{A,B,C,D\}$. For $k\in\NM$ consider the periodic configurations $\eta^x_k:=S^k(x^\infty)$ with corresponding periodic subshifts $\Xi^x_k:=\Orb(\eta^x_k)$. Then $\Xi$ is periodically approximable and the sequences of periodic subshifts $(\Xi^x_k)_k$ converge to $\Xi$.
\end{proposi}
%%%%%%%%%%%%%%%%

\noindent {\bf Proof: } Consider the closed paths

$$\gamma_A:= (AB,BA)\,,\quad
	\gamma_B:= (AC,CA)\,,\quad
	\gamma_C:= (DB,BD)\,,\quad
	\gamma_D:= (DC,DC)\,,
$$

\noindent in the GAP-graph $\gs_1$ of order $1$, c.f. Figure~\ref{1D.fig-GAPgrs}. The associated periodic word (Definition~\ref{1D.def-AssPerWor}) satisfy

$$\eta(\gamma_A):= S(A^\infty)\,,\quad
	\eta(\gamma_B):= S(B^\infty)\,,\quad
	\eta(\gamma_C):= S(C^\infty)\,,\quad
	\eta(\gamma_D):= S(D^\infty)\,.
$$

\noindent Thus, the statement follows from Proposition~\ref{1D.prop-PrimSub} since the substitution $S$ is primitive.
\hfill $\Box$

\vspace{.2cm}

\noindent Combining Proposition~\ref{1D.prop-GRS} with Theorem~\ref{1D.th-ContSpectrGenSchrDyn} leads to

%%%%%%%%%%%%%%%%
\begin{coro}
\label{1D.cor-GRS}
Let $\Xi$ be the Golay-Rudin-Shapiro subshift and $H$ be a Hamiltonian defined by Equation~\ref{1D.eq-schdis} satisfying (R1, R2, R3). Then the equations

$$\sigma(H_\xi)
	=\lim_{k\to\infty} \sigma(H_{\eta^x_k})
$$

\noindent hold for every $\xi\in\Xi$ and $x\in\{A,B,C,D\}$.
\end{coro}
%%%%%%%%%%%%%%%%

%%%%%%%%%%%%%%%%%%%%%%%%%%%%%%%%%%%%%%%%%%%%%%%%%%%%%%%%%%%%%%%%%%%%%%%%%
%%%%%%%%%%%%%%%%%%%%%%%%%%%%%%%%%%%%%%%%%%%%%%%%%%%%%%%%%%%%%%%%%%%%%%%%%
\appendix

\section{Pattern equivariant algebra}
\label{1D.sect-PattEqOp}

\noindent Pattern equivariant functions were defined in \cite{Ke03}. In this appendix, the pattern equivariant algebra will be defined by finite range operators on $\ell^2(\ZM)$ with pattern equivariant coefficients. 
Such an algebra plays, in the noncommutative approach advocated in \cite{BBdN17}, a role analog to the algebra of polynomials on an algebraic variety. The rigidity of these functions with respect to the local structures makes the analysis of the related operators more difficult.

%%%%%%%%%%%%%%%%%%
\begin{defini}[Pattern equivariant function \cite{Ke03}]
\label{1D.def-PattEqFunct}
A function $g:\csp\to\CM$ is called {\em (strongly) pattern equivariant} if there exists a radius $r\in\NM$ such that $g(\xi)=g(\eta)$ for all $\xi,\eta\in\csp$ with $\xi|_{[-r,r-1]} = \eta|_{[-r,r-1]}$.
\end{defini}
%%%%%%%%%%%%%%%%%%

\noindent For the following assertion, recall the notion of $\os(u,v)\subseteq\csp$ introduced Equation~\eqref{1D.eq-clopen}. The set $\Cc(\csp)$ of continuous functions on $\csp$ defines a \Cs if equipped with the pointwise multiplication and the uniform norm $\|\cdot\|_\infty$.

%%%%%%%%%%%%%%%%%%
\begin{proposi}
\label{1D.prop-CharPatEq}
For $g:\csp\to\CM$ the following assertions are equivalent.
\begin{itemize}
\item[(a)] The function $g$ is (strongly) pattern equivariant.
\item[(b)] There is an $N\in\NM$, coefficients $c_j\in\CM$ and $u_j,v_j\in\as^\ast$ for $1\leq j\leq N$ such that

$$g
	=\sum_{j=1}^N c_j\;\chi_{\os(u_j,v_j)}.
$$
\item[(c)] The function $g:\csp\to\CM$ is continuous and takes finitely many values.
\end{itemize}
In particular, the set of (strongly) pattern equivariant functions $\pea\subseteq\Cc(\csp)$ are a dense $\ast$-subalgebra.
\end{proposi}
%%%%%%%%%%%%%%%%%%

\noindent {\bf Proof: } Sets of the form $\os(u,v)\subseteq\csp$ for $u,v\in\as^\ast$ define a base of the product topology on $\csp$, c.f. Equation~\eqref{1D.eq-clopen}. Thus, the set of functions of the form in (b) are clearly dense in $\Cc(\csp)$ with respect to the uniform norm. Moreover, it is immediate to check that the pointwise product $\chi_{\os(\tu,\tv)}\cdot\chi_{\os(u',v')}$ equals either $0$ or $\chi_{\os(u,v)}$ for some $u,v\in\as^\ast$. Hence, these functions form a $\ast$-subalgebra of $\Cc(\csp)$. It is left to show that (a), (b) and (c) are equivalent:

\vspace{.1cm}

\noindent (a)$\Rightarrow$(b): Let $g:\csp\to\CM$ be (strongly) pattern equivariant. Then there is an $r\in\NM$ with $g(\xi)=g(\eta)$ for all $\xi,\eta\in\csp$ satisfying $\xi|_{[-r,r-1]} = \eta|_{[-r,r-1]}$. Consider the finite set 

$$\{(u,v) \;:\; uv\in\as^\ast,\; |u|=|v|=r\}
	=\{(u_1,v_1),\ldots,(u_N,v_N)\}\,.
$$ 

\noindent Define $c_j:=g(\xi)$ for one $\xi\in\os(u_j,v_j)$. As $g$ is pattern equivariant $c_j$ is independent of the choice of $\xi$. Since $\csp$ equals to the disjoint union of the sets $\os(u_j,v_j)\,,\; 1\leq j\leq N\,,$ we get

$$g
	=\sum_{j=1}^N c_j\;\chi_{\os(u_j,v_j)}
$$

\noindent (b)$\Rightarrow$(c): This is clear as the sets of the form $\os(u_j,v_j)$ are clopen.

\vspace{.1cm}

\noindent (c)$\Rightarrow$(a): Let $g:\csp\to\CM$ be continuous and suppose it takes only finitely many values denoted by $c_1,\ldots,c_N$. Let $U_j\subseteq\CM$ be an open neighborhood of $c_j$ such that $U_j\cap U_i=\emptyset$ if $i\neq j$. By continuity, the preimage $g^{-1}(U_j)=g^{-1}(\{c_j\})\subseteq\csp$ is clopen and these sets form a finite partition of $\csp$. Consequently, if $\xi\in\csp$ satisfies $g(\xi)=c_j$, there is an $r'\in\NM$ such that if $\eta\in\os(\xi|_{[-r',-1]},\xi|_{[0,r'-1]})$, then $g(\eta)=g(\xi)=c_j$. Since such an $r'\in\NM$ exists for every $\xi\in g^{-1}(U_j)$ and $\os(\xi|_{[-r',-1]},\xi|_{[0,r'-1]})$ is open, there is a maximal $r_j\in\NM$ satisfying $g(\xi)=g(\eta)=c_j$ whenever $\xi|_{[-r_j,r_j-1]} = \eta|_{[-r_j,r_j-1]}$ by using the compactness of $g^{-1}(U_j)$. Hence, $g$ is strongly pattern equivariant with $r:=\max_{1\leq j\leq N} r_j$.
\hfill$\Box$

\vspace{.2cm}

\noindent Clearly, the statement of Proposition~\ref{1D.prop-CharPatEq} is valid if $\csp$ is replaced by any $\Xi\in\Inv$. We denote by $\pea(\Xi)$ the set of all (strongly) pattern equivariant functions $g:\Xi\to\CM$.

\vspace{.1cm}

\noindent The shift $\tra$, acting on any subshift $\Xi$, is an homeomorphism of $\Xi$. It generates a $\ZM$ action on $\Xi$. In particular it defines the groupoid $\Xi\rtimes_{\tri}\ZM$, called the crossed product (for a comprehensive review and details, see \cite{Re80, Bec16,BBdN17}). By duality, the shift defines an action on the \Cs $\Cc(\Xi)$ of continuous functions on $\Xi$, defined by $\phi(f)=f\circ\tra^{-1}$. Clearly the space $\pea(\Xi)$ is a dense $\ast$-subalgebra of $\Cc(\Xi)$. It is easy to check that $\pea(\Xi)$ is $\phi$-invariant. By construction, the algebra $\pea(\Xi)\rtimes_{alg}\ZM$, is the algebra generated by $\pea(\Xi)$ and a unitary element $S$ such that 

$$S\,f\,S^{-1}=f\circ\tra\,.
$$

\noindent Using the left regular representation on $\ell^2(\ZM)$, the reader is invited to check that this algebra is a $\ast$-algebra made of operators of the form provided in Equation~\ref{1D.eq-schdis} (see also \cite{Bec16}, Theorem~3.7.10). Namely they have finite range and pattern equivariant coefficients. Using the construction described in \cite{Re80, Bec16,BBdN17}, this algebra is a dense $\ast$-subalgebra of the \Cs $\Cc_{red}^\ast(\Xi\rtimes_{\tra}\ZM)$. It is important to remark at this point, that, since $\ZM$ is an amenable group, $\Cc(\Xi)\rtimes\ZM$ coincides with \Cs $\Cc^\ast_{red}(\Xi\rtimes_{\tra}\ZM)$. This can be summarized as

%%%%%%%%%%%%%%%%
\begin{coro}
\label{1D.cor-PEden}
Let $\Xi\in\Inv$. Then $\pea(\Xi)\rtimes_{alg}\ZM$ is a dense $\ast$-subalgebra of $\Cc(\Xi)\rtimes\ZM$.
\end{coro}
%%%%%%%%%%%%%%%%

%%%%%%%%%%%%%%%%%%%%%%%%%%%%%%%%%%%%%%%%%%%%%%%%%%%%%%%%%%%%%%%%%%%%%%%%%

%%%%%%%%%%%%%%%%%%%%%%%%%%% September 27, 2016 %%%%%%%%%%%%%%%%%%%%%%%%%%%
%%%%                                                                  %%%%
%%%%       Spectral continuity for aperiodic quantum systems II.      %%%%
%%%%                  Periodic approximations in 1D                   %%%%
%%%%                                                                  %%%%
%%%%      Siegfried Beckus, Jean Bellissard, Giuseppe De Nittis,      %%%%
%%%%                   Revised Version #6: 02/18/18                   %%%%
%%%%                                                                  %%%%
%%%%%%%%%%%%%%%%%%%%%%%%%%%%%%%%%%%%%%%%%%%%%%%%%%%%%%%%%%%%%%%%%%%%%%%%%%

%%%%%%%%%%%%%%%%%%%%%%%%%%%%%%%%%%%%%%%%%%%%%%%%%%%%%%%%%%%%%%%%%%%%%%%%%


\begin{thebibliography}{9999}

%\bibitem{}

\bibitem{AS93} J.-P.~Allouche, J.~O.~Shallit, ``Complexit\'e des suites de Rudin-Shapiro g\'en\'eralis\'ees'', {\em . Th\'eorie Nombres, Bordeaux}, {\bf 5}, 283-302, (1993).

\bibitem{Al94} J.-P.~Allouche, ``Sur la complexit\'e des suites infinies, {\em Bull.  Belg.  Math.  Soc.},  {\bf 1}, 133-143, (1994).

\bibitem{Al97} J.-P.~Allouche, ``Schr\"odinger operators with Rudin-Shapiro potentials are not palindromic. Quantum problems in condensed matter physics'', {\em J. Math. Phys.}, {\bf 38}, 1843–1848, (1997).

\bibitem{AP98} J.~E.~Anderson, I.~F.~Putnam, ``Topological invariants for substitution tilings and their associated $C^\ast$-algebra'', in {\em Ergod. Th, \& Dynam. Sys.}, {\bf 18}, 509-537, (1998).

\bibitem{AR91} P.~Arnoux, G.~Rauzy, ``Repr\'esentation g\'eom\'etrique de suites de complexit\'e $2n+1$, {\em Bull.  Soc. Math. France}, {\bf 199}, 199-215, (1991).

\bibitem{AG95} F.~Axel, D.~Gratias, {\em Beyond Quasicrystals}, Les Houches, March 7–18, 1994, Springer-Verlag Berlin Heidelberg, (1995).

\bibitem{Bec16} S.~Beckus, ``Spectral approximation of aperiodic Schr\"odinger operators'', PhD Thesis, Friedrich-Schiller-Universit\"at Jena, (2016).

\bibitem{BB16} S.~Beckus, J.~Bellissard, ``Continuity of the Spectrum of a Field of Self-Adjoint Operators'', {\em Annales Henri Poincar{\'e}}, {\bf 17}, 3425-3442, (2016).

\bibitem{BBdN17} S.~Beckus, J.~Bellissard, G.~de~Nittis ``{S}pectral continuity for aperiodic quantum systems {I}. {G}eneral theory'', arXiv:1709.00975, (2017).

\bibitem{BBdN18} S.~Beckus, J.~Bellissard, G.~de~Nittis ``{S}pectral continuity for aperiodic quantum systems: Book of examples'', preprint, (2018).

\bibitem{BP17} S.~Beckus, F.~Pogorzelski ``{D}elone dynamical systems and spectral convergence'', arXiv:1711.07644, (2017).

\bibitem{Be86} J.~Bellissard, ``$K$-Theory of $C^{\ast}$-algebras in Solid State Physics'', in {\em Statistical Mechanics and Field Theory,  Mathematical Aspects}, T.C. Dorlas, M.~N.~Hugenholtz \& M.~Winnink Eds., Lecture Notes in Physics, {\bf 257}, 99-156, (1986).

\bibitem{BIST89} J.~Bellissard, B.~Iochum, E.~Scoppola, D.~Testard, ``Spectral properties of one-dimensional quasi-crystals'', in {\em Comm. Math. Phys.}, {\bf 125}, no. 3, 527-543, (1989).

\bibitem{Be90} J.~Bellissard, ``Spectral properties of Schr\"odinger's operator with a Thue-Morse potential'', in {\em Number theory and physics} (Les Houches, 1989), 140–150, Springer Proc. Phys., {\bf 47}, Springer, Berlin, 1990.

\bibitem{BBG91} J.~Bellissard, A.~Bovier, J.-M.~Ghez, ``Spectral properties of a tight binding Hamiltonian with period doubling potential'', {\em Comm. Math. Phys.}, {\bf 135}, 379-399, (1991).

\bibitem{BIT91} J.~Bellissard, B.~Iochum, D.~Testard, ``Continuity properties of the electronic spectrum of 1D quasicrystals'', {\em Commun. Math. Phys.}, {\bf 141}, 353–380, (1991).

\bibitem{Bel03} J.~Bellissard, ``Coherent and dissipative transport in aperiodic solids'', Published in  {\em Dynamics of Dissipation}, P. Garbaczewski, R. Olkiewicz (Eds.), Lecture Notes in Physics, {\bf 597}, Springer (2003), pp.  413-486, (2003), (Proceedings of the 38th Winter School of Theoretical Physics, Ladek, Poland, 6-15 Feb 20023.

\bibitem{Bel13} J.~Bellissard, ``Bloch Theory for 1-D FLC Aperiodic Media'', Lecture given at WCAOS, UC Davis, October 26th, 2013. See slides at {\verb+http://people.math.gatech.edu/~jeanbel/TalksE/wannier13.pdf+}

\bibitem{Be15} J.~Bellissard, ``Delone Sets and Material Science: a Program'', in  {\em Mathematics of Aperiodic Order}, J. Kellendonk, D. Lenz, J. Savinien, Eds., "Progress in Mathematics" series, 309, Birkh\"auser, 2015.

\bibitem{BS91} V.~G.~Benza, C.~Sire, ``Electronic  spectrum of the octagonal quasicrystal: chaos, gaps and level clustering'', {\em Phys. Rev. B}, {\bf 44}, 10343-10345, (1991).

\bibitem{BK13} G.~Berkolaiko,P.~Kuchment, {\em Introduction to Quantum Graphs}, American Mathematical Society, Mathematical Survey and Monographs, vol. {\bf 186}, (2013).

\bibitem{BG93} A.~Bovier, J.-M.~Ghez, ``Spectral properties of one-dimensional Schr\"odinger operators with potentials generated by substitutions'', {\em Commun. Math. Phys.}, {\bf 158}, 45-66, (1993).
\bibitem{Bru46} N.~G.~de~Bruijn, ``A combinatorial problem'', {\em Nederl. Akad. Wetensch}, {\bf 49}, 758-764, (1946).

\bibitem{CKM87} R.~Carmona, A.~Klein, F.~Martinelli, ``Anderson localization for Bernoulli and other singular potentials'',
 {\em Commun. Math. Phys.}, {\bf 108}, 41-66, (1987).

\bibitem{Ca86} M.~Casdagli, ``Symbolic dynamics for the renormalization map of a quasiperiodic Schr\"odinger equation'', {\em Comm. Math. Phys.}, {\bf 107}, 295-318, (1986).

\bibitem{Ca97} J.~Cassaigne, ``Complexit\'e et facteurs sp\'eciaux'', {\em Bull. Belg. Math. Soc. Simon Stevin}, (\textbf{4}), Springer-Verlag, no. 1, 67--88, (1997).

\bibitem{Ch50} C.~Chabauty,``Limite d'ensembles et g\'eom\'etrie des nombres'', in {\em Bull. Soc. Math. France}, {\bf 78}, 143-151, (1950).

\bibitem{Ce37} E.~\v{C}ech,``On bicompact spaces'', {\em Ann. of Math.}, {\bf 2}, 823-844, (1937).

\bibitem{CH73} E.~M.~Coven, G.~A.~Hedlund, ``Sequences with minimal block growth'', {\em Math. Systems Theory}, {\bf 7}, 138-153, (1973).

\bibitem{Da98} D.~Damanik, ``Singular continuous spectrum for the period doubling {H}amiltonian on a set of full measure'', {\em Comm. Math. Phys.}, {\bf 196}, no. 2, 477-483, (1998).

\bibitem{Da00} D.~Damanik, ``Local symmetries in the period-doubling sequence'', {\em Discrete Appl. Math.}, {\bf 100}, 115-121, (2000).

\bibitem{Da00b} D.~Damanik, ``Substitution {H}amiltonians with bounded trace map orbits'', {\em J. Math. Anal. Appl.}, {\bf 249}, 393-411, (2000).

\bibitem{DL03} D.~Damanik, D.~Lenz, ``Uniform spectral properties of one-dimensional quasicrystals. {IV}. {Q}uasi-{S}turmian potentials'', {\em J. Anal. Math.}, {\bf 90}, 115-139, (2003).

\bibitem{DL06} D.~Damanik, D.~Lenz, ``Zero-measure {C}antor spectrum for {S}chr\"odinger operators with low-complexity potentials'', {\em J. Math. Pures Appl. (9)}, {\bf 85}, 671-686, (2006).

\bibitem{DLQ14} D.~Damanik, Q.-H.~Liu, Y.-H.~Qu, ``Schr\"odinger Operators with Dynamically Defined Potentials: A Survey'', {\tt arXiv:1410.2445}, (2014).

\bibitem{DEGT08} D.~Damanik, M.~Embree, A.~Gorodetski, S.~Tcheremchantsev, ``The fractal dimension of the spectrum of the {F}ibonacci {H}amiltonian'', {\em Comm. Math. Phys.}, {\bf 280}, 499-516, (2008).

\bibitem{DEG15} D.~Damanik, M.~Embree, A.~Gorodetski, ``Spectral properties of Schr\"odinger operators arising in the study of quasicrystals'', in {\em Mathematics of aperiodic order}, J. Kellendonk, D. Lenz, J. Savinien, Eds., "Progress in Mathematics" series, 309, Birkh\"auser, 2015.

\bibitem{DGS15} D.~Damanik, A.~Gorodetski, B.~Solomyak, ``Absolutely continuous convolutions of singular measures and an application to the square Fibonacci Hamiltonian'', {\em Duke Math. J.}, {\bf 164}, 1603-1640, (2015).

\bibitem{DGY16} D.~Damanik, A.~Gorodetski, W.~Yessen, ``The Fibonacci Hamiltonian'', {\tt arXiv:1403.7823}, March 2014. to appear in {\em Inventiones Mathematicae}.

\bibitem{DD63} J.~Dixmier, A.~Douady, ``Champs continus d'espaces hilbertiens et de $C^\ast$-alg\`ebres'', {Bull. Soc. Math. France}, {\bf 91}, 227-284, (1963).

\bibitem{Di69} J.~Dixmier, {\em Les $C^{\ast}$-alg\`ebres et leurs repr\'esentations}, (French) Deuxi\`eme \'edition. Cahiers Scientifiques, Fasc. XXIX. Gauthier-Villars \'Editeur, Paris 1969; (ii) J.~Dixmier, {\em C$^\ast$-algebras}, North-Holland, Amsterdam-New York-Oxford, (1977).

\bibitem{ER82} A.~Ehrenfeucht, G.~Rozenberg, ``On subword complexities of homomorphic images of languages'', {\em RAIRO Inform. Th\'{e}or.}, {\bf 16}, 303-316, (1982).

\bibitem{Fe62} J.~M.~G.~Fell, ``A Hausdorff topology for the closed subsets of a locally compact non-Hausdorff space'', {\em Proc. Amer. Math. Soc.}, {\bf 13}, 472-476, (1962).

\bibitem{Fi14} F.~Fiorenzi, ``Periodic configurations of subshifts on groups'', arXiv:1402.3448, (2014).

\bibitem{Fl94} C.~Flye Sainte Marie, ``Question 48'', {\em L'Interm\'ediare Math.}, {\bf 1}, 107-110, (1894).

\bibitem{Fog02} P.~Fogg, ``Substitutions in dynamics, arithmetics and combinatorics'', Lecture Notes in Mathematics, Springer, (2002).

\bibitem{Ga02} F.~G{\"a}hler, talk given at the Conference {\em Aperiodic Order, Dynamical Systems, Operator Algebra and Topology}, August 4-8, 2002, Victoria, B.C., Canada. {\em (unpublished)}.

\bibitem{Go51} M.~J.~L.~Golay, ``Static multilist spectrometry and its application to panoramic display of infrared spectra'', {\em J. Optical Soc. America}, {\bf 41}, 468-472, (1951).

\bibitem{Goo46} I.~J.~Good, ``Normal recurring decimals'', {\em J. London Math Soc.}, {\bf 21}, 167-169, (1946).

\bibitem{GH55} W.~H.~Gottschalk, G.~A.~Hedlund, {\em Topological dynamics}, American Mathematical Society, Colloquium Publications, Vol. 36, Providence, Rhode Island, 1955.

\bibitem{Gr69} F.~P.~Greenleaf, {\em Invariant Means on Topological Groups and Their Applications}, Van Nostrand Reinhold, (1969).

\bibitem{HM40} G.~A.~Hedlund, M.~Morse, ``Symbolic Dynamics II: Sturmian trajectories'', {\em Amer. J. Math.}, {\bf 62}, 1-42, (1940).

\bibitem{HKS95} A.~Hof, O.~Knill, B.~Simon, ``Singular continuous spectrum for palindromic Schrödinger operators'', {\em Commun. Math. Phys}, {\bf 174}, 149-159, (1995).

\bibitem{Ju10} A.~Julien, ``Complexity and cohomology for cut-and-projection tilings'', {\em Ergodic Theory and Dynamical Systems}, {\bf 30}, 489-523, (2010).

\bibitem{Ke03} J.~Kellendonk, ``Pattern-equivariant functions and cohomology'', {\em J. Phys. A}, {\bf 36}, 5765-5772, (2003).

\bibitem{KP17} J.~Kellendonk, E.~Prodan, ``Bulk-boundary correspondance for Sturmian Kohmoto like models'',\\ {\tt arXiv:1710.07681}, (2017).

\bibitem{KKT83} M.~Kohmoto, L.~Kadanoff, C.~Tang, ``Localization problem in one dimension: Mapping and escape'', {\em Phys. Rev. Lett.}, {\bf 50}, 1870–1872, (1983).

\bibitem{KO84} M.~Kohmoto, Y.~Oono, ``Cantor spectrum for an almost periodic Schr\"odinger equation and a dynamical map'', {\em Phys. Letters A}, {\bf 102}, 145-148, (1984).

\bibitem{KO89} S.~Kotani, ``Jacobi matrices with random potentials taking finitely many values, {\em Rev. Math. Phys.}, {\bf 1}, 129-133, (1989).

\bibitem{La99A} J.~C.~Lagarias, ``Geometric models for quasicrystals I. Delone sets of finite type''. {\it Discrete Comput. Geom.}, {\bf 21}, 161-191, (1999).

\bibitem{Le02} D.~Lenz, ``Singular spectrum of {L}ebesgue measure zero for one-dimensional quasicrystals''. {\it Comm. Math. Phys.}, {\bf 227}, 119-130, (2002).

\bibitem{LM95} D.~Lind, B.~Marcus, ``An introduction to symbolic dynamics and coding'', Cambridge University Press, Cambridge, (1995).

\bibitem{LQ11} Q.~H.~Liu, Y.~H.~Qu ``Uniform convergence of {S}chr\"odinger cocycles over simple {T}oeplitz subshift'', {\em Ann. Henri Poincar\'e}, {\bf 12}, 153-172, (2011).

\bibitem{LQ15} Q.~H.~Liu, Y.~H.~Qu ``On the {H}ausdorff dimension of the spectrum of the {T}hue-{M}orse {H}amiltonian''. {\it Comm. Math. Phys.}, {\bf 338}, 867-891, (2015).

\bibitem{LQY16} Q.~H.~Liu, Y.~H.~Qu, X.~Yiao, ``''Mixed spectral nature'' of the Thue-Morse Hamiltonian'',\\ {\tt arXiv:1512.08011}, (2016).

\bibitem{Mo05} E.~Moreno, ``De Bruijn sequences and de Bruijn graphs for a general language'', {\em . Inform. Process. Lett. 96}, {\bf 6}, 214-219, (2005).

\bibitem{OPRSS83} S.~Ostlund, R.~Pandit, D.~Rand, H.~J.~Schellnhuber, E.~D.~Siggia, ``One-Dimensional Schr\"odinger Equation with an Almost Periodic Potential'', {\em Phys. Rev. Lett.}, {\bf 50}, 1873-1876, (1983).

\bibitem{OK85} S.~Ostlund, S.-H. Kim, ``Renormalization of Quasiperiodic Mappings'', {\em Physica Scripta}, {\bf T9}, 193-198, (1985).

\bibitem{Pr13} E.~Prodan, ``Quantum transport in disordered systems under magnetic fields: A study based on operator algebras'', {\em Appl. Math. Res Express}, {\bf Vol. 2013}, 176-255, (2013).

\bibitem{Qu87} M.~Queff{\'e}lec, ``Substitution dynamical systems--spectral analysis'', Lecture Notes in Math. {\bf 1294}, Berlin, Springer-Verlag (1987).

\bibitem{Qu10} M.~Queff{\'e}lec, ``Substitution dynamical systems--spectral analysis'', Lecture Notes in Math. {\bf 1294}, 2nd Edition, Berlin, Springer-Verlag (2010).

\bibitem{RW92} C.~Radin, M.~Wolff, ``Space tilings and local isomorphism'', {\em Geom. Dedicata}, {\bf 42}, no. 3, Edited by V. Berthי, S. Ferenczi, C. Mauduit and A. Siegel, 355--360, (1992).

\bibitem{Rau83} G.~Rauzy, `` Suites \`a termes dans un alphabet fini'', Seminar on number theory, 1982-1983, Exp. No. 25, Univ. Bordeaux I, Talence (1983).

\bibitem{Re80} J.~Renault, ``A Groupoid Approach to $C^{\ast}$-Algebras'', Lecture Notes in Math., {\bf 793}, Springer, Berlin, (1980).

\bibitem{Rud59} W.~Rudin, ``Some theorems on Fourier coefficients'', {\em Proc. Amer. Math. Soc.}, {\bf 10}, 855-859, (1959).

\bibitem{Sa03} L.~Sadun, ``Tiling spaces are inverse limits'', {\em J. Math. Phys.}, {\bf 44}, no. 11, 5410-5414, (2003).

\bibitem{Sa08} L.~Sadun, ``Topology of tiling spaces'', {\em American Mathematical Society, Providence, RI}, University Lecture Series, {\bf 46}, (2008).

\bibitem{Sh51} H.~S.~Shapiro, ``Extremal problems for polynomials and power series'', Master's thesis, M.I.T., Cambridge Mass., 1951.

\bibitem{SBGC84} D.~Shechtman, I.~Blech, D.~Gratias, J.~W.~Cahn, ``Metallic Phase with Long-Range Orientational Order and No Translational Symmetry'', {\em Phys. Rev. Lett.}, {\bf 53}, 1951-1953, (1984).

\bibitem{Su87} A.~S\"ut\H{o}, ``The spectrum of a quasi-periodic Schr\"odinger operator'', {\em Commun. Math. Phys.}, {\bf 111}, 409-415, (1987).

\bibitem{Su89} A.~S\"ut\H{o}, ``Singular continuous spectrum on a {C}antor set of zero {L}ebesgue measure for the {F}ibonacci {H}amiltonian'' , {\em J. Statist. Phys.}, {\bf 56}, 525-531, (1989).

\bibitem{Su95} A.~S\"ut\H{o}, ``Schr\"odinger difference equation with deterministic ergodic potentials'', in {\em Beyond Quasicrystals}, Les Houches, March 7–18, 1994, Springer-Verlag Berlin Heidelberg, (1995).

\bibitem{Ty30} A.,~Tychonoff, ``\"{U}ber die topologische {E}rweiterung von {R}\"aumen'', {\em Math. Ann.}, {\bf 102}, no. 1, 544-561, (1930).

\bibitem{Ur24} P.,~Urysohn, ``\"{U}ber die {M}etrisation der kompakten topologischen {R}\"aume'', {\em Math. Ann.}, { \bf 92}, no. 3-4, 275-293, (1924).

\bibitem{Vi22} L.~Vietoris, ``Bereiche zweiter {O}rdnung'', in {\em Monatsh. Math. Phys.}, {\bf 32}, no. 1, 258-280, (1922).

\bibitem{Que01}  B.~von~Querenburg, ``Mengentheoretische Topologie'', 3. Auflage, Springer-Verlag, Berlin-New York (2001).

\bibitem{Wa82} P.~Walters, ``An introduction to ergoodic theory'',  {\bf 79}, Springer-Verlag (1982).

\end{thebibliography}
\end{document}